\pdfoutput=1
\documentclass[11pt,a4paper]{amsart}

\usepackage[a4paper,margin=1in]{geometry}
\usepackage{amsmath,amssymb,amsthm,mathtools}
\usepackage{graphicx}
\usepackage{booktabs}
\usepackage{microtype}
\usepackage{xcolor}
\usepackage[font=small,labelfont=bf]{caption}
\usepackage[numbers,sort&compress]{natbib}
\bibpunct[, ]{[}{]}{,}{n}{,}{,}
\usepackage[colorlinks=true,linkcolor=blue!55!black,citecolor=blue!55!black,
  urlcolor=blue!55!black,
  pdftitle={A One-Sided Chord Trim of Wetzel's Sector Covers Unit Arcs},
  pdfauthor={Chatchawan Panraksa and Wacharin Wichiramala}]{hyperref}

\theoremstyle{plain}
\newtheorem{theorem}{Theorem}[section]
\newtheorem{lemma}[theorem]{Lemma}
\newtheorem{proposition}[theorem]{Proposition}
\newtheorem{corollary}[theorem]{Corollary}
\newtheorem*{wedgelemma}{Wedge Lemma}
\theoremstyle{definition}
\newtheorem{definition}[theorem]{Definition}
\newtheorem{example}[theorem]{Example}
\newtheorem{openproblem}[theorem]{Open Problem}
\theoremstyle{remark}
\newtheorem{remark}[theorem]{Remark}

\newcommand{\R}{\mathbb{R}}
\newcommand{\conv}{\operatorname{conv}}
\newcommand{\area}{\operatorname{area}}
\newcommand{\ip}[2]{\langle #1,#2\rangle}
\newcommand{\pt}[1]{p_{#1}}
\newcommand{\PiSector}{\Pi}
\newcommand{\refl}{\sigma}
\newcommand{\sC}{\mathsf{C}}
\newcommand{\sR}{\mathsf{R}}
\newcommand{\arcZ}{\mathcal{Z}}

\begin{document}

\title[A one-sided chord trim of Wetzel's sector]{A One-Sided Chord Trim of Wetzel's Sector Covers Unit Arcs}

\author{Chatchawan Panraksa}
\address{Applied Mathematics Program, Science Division, Mahidol University
International College, 999 Phutthamonthon 4 Road, Salaya, Nakhon Pathom
73170, Thailand}
\email{chatchawan.pan@mahidol.ac.th}
\urladdr{https://orcid.org/0000-0003-0692-5453}

\author{Wacharin Wichiramala}
\address{Department of Mathematics and Computer Science, Faculty of Science,
Chulalongkorn University, 254 Phayathai Road, Pathumwan, Bangkok 10330,
Thailand}
\email{wacharin.w@chula.ac.th}

\date{November 12, 2026}

\begin{abstract}
Wetzel conjectured, and Panraksa and Wichiramala proved, that the
circular sector of radius one and angle $30^\circ$ contains a congruent
copy, reflections allowed, of every arc of length one. We show that half
of its arc is not needed: replacing the arc between the angles $15^\circ$
and $30^\circ$ by its chord gives a convex cover of area
$\pi/24+(\sqrt6-\sqrt2)/8=0.2603092\ldots$. As in the proof for the
sector, we place an arc in two supporting wedges. A sum of support, cap
and radial inequalities, with vector coefficients summing to zero, shows
that an arc with no copy in the trimmed sector is longer than one.
One two-parameter family of configurations is checked by rational
interval arithmetic; programs and data are supplied. No chord that
starts below $3.6^\circ$ gives a cover. We also reduce the
inclusion-minimality of the trimmed sector to two families of cuts.
\end{abstract}

\maketitle

\section{Introduction}
\label{sec:intro}

A \emph{unit arc} is the image of a continuous rectifiable map
$\gamma:[0,1]\to\R^2$ of total length one. A set $K\subset\R^2$ is a
\emph{cover} if it contains a congruent copy of every unit arc; as in
\cite{wetzel-wichiramala2010,panraksa-wichiramala2021}, congruent copies
may be reflected. Moser's worm problem asks for the least area of a cover
\cite{moser1966,moser1991}; see \cite[Problem~D18]{croft-falconer-guy1991}
and \cite[Section~11.4]{brass-moser-pach2005}. The problem is open even
for convex covers; Wetzel \cite{wetzel2003} and Movshovich
\cite{movshovich2025} survey its history.

Wetzel \cite{wetzel1973} studied circular sectors that cover all curves
of a given length. He conjectured that the sector of radius one and angle
$30^\circ$ is a cover
\cite{gerriets-poole1974,norwood-poole-laidacker1992,wetzel2003}, and we
proved this in \cite{panraksa-wichiramala2021}. Its area is
$\pi/12\approx0.2617994$. We ask whether a proper convex subset of the
sector is still a cover.

The question goes back to Wetzel, who predicted that a set can be clipped
from the sector at one end of its circular arc so as to leave a cover,
which is then no longer symmetric; see Movshovich's survey
\cite[Section~4]{movshovich2025}. The survey also conjectures the shape
of a convex cover of least area: a subset of the sector bounded by one
radius, a segment of the other radius that starts at the apex, and a
convex arc that joins the end of this segment to the end of the first
radius and meets the circular arc of the sector nowhere else; the two
ends are conjectured to be at distance $1/2$. Both statements are about
covers of all unit arcs. As far as we know, neither has been proved, and
no proper subset of the sector has been shown to be a cover.

We keep the circular boundary from angle $0^\circ$ to $15^\circ$ and
replace the rest of the arc by its chord. Write $O$ for the origin and
$\pt\theta=(\cos\theta,\sin\theta)$, with angles in degrees, and put
\[
 \PiSector=\{r\pt\theta:0\le r\le1,\ 0^\circ\le\theta\le30^\circ\},
 \qquad
 H_{15}=\PiSector\cap
 \{x:\ip{x}{\pt{22.5}}\le\cos7.5^\circ\}.
\]
The boundary of $H_{15}$ consists of the radii $[O,\pt0]$ and
$[O,\pt{30}]$, the circular arc from $\pt0$ to $\pt{15}$, and the chord
$[\pt{15},\pt{30}]$ (Figure~\ref{fig:trim}(a)). We call $H_{15}$ the
\emph{trimmed sector}.

\begin{theorem}
\label{thm:H15}
Every unit arc has a congruent copy contained in $H_{15}$. The set
$H_{15}$ is compact and convex, and
\[
 \area(H_{15})=\frac{\pi}{24}+\frac{\sqrt6-\sqrt2}{8}=0.2603092\ldots.
\]
\end{theorem}

The set conjectured in \cite{movshovich2025} omits much more of the
boundary of the sector than $H_{15}$ does. Placed so that the radius it
keeps is $[O,\pt0]$, it meets the circular arc only at $\pt0$ and does
not contain $\pt{30}$, whereas $H_{15}$ keeps both radii and the arc from
$\pt0$ to $\pt{15}$. We return to that set in
Section~\ref{subsec:minimality}.

By Theorem~\ref{thm:H15}, the infimum of the areas of convex covers is at
most $\area(H_{15})<0.2603093$; the best published lower bound is
$0.232239$, due to Khandhawit, Pagonakis and Sriswasdi
\cite{khandhawit-pagonakis-sriswasdi2013}, which improves the bound of
Khandhawit and Sriswasdi \cite{khandhawit-sriswasdi2007}. The area of
$H_{15}$ is $0.001490$ below that of the sector and $0.000647$ below that
of the scaled $30^\circ$--$60^\circ$--$90^\circ$ triangle of our
computer-assisted preprint
\cite[Proposition~20]{wichiramala-panraksa2026}. Earlier
convex covers include one of area less than $0.28610$
\cite[Theorem~2]{gerriets-poole1974} and one of area about $0.275237$
\cite{norwood-poole-laidacker1992}. The nonconvex covers of Norwood and
Poole \cite{norwood-poole2003}, who allow only translations and rotations,
and of Ploymaklam and Wichiramala \cite{ploymaklam-wichiramala2018} have
areas $0.260437$ and $0.26007$, respectively. Two recent preprints
announce smaller convex covers of other shapes: a triangle of area
$0.2578835951\ldots$ \cite[Corollary~7.2]{deng2026} and a quadrilateral
of area $0.2463322372\ldots$, together with the lower bound $0.239$ for
the area of a convex cover \cite[Theorem~1.1]{keller2026}. Neither fits
inside $H_{15}$, which has diameter one, since the two equal sides of
Deng's triangle and a diagonal of Keller's quadrilateral are longer than
one \cite[Section~7]{deng2026}, \cite[Section~4]{keller2026}.
Table~\ref{tab:areas} lists Wetzel's sector and the smaller convex covers
mentioned above. Theorem~\ref{thm:H15} is therefore not an area record.
Its interest lies in the fact that a proper convex subset of Wetzel's
sector is still a cover, and in the method of proof.

Among convex covers of diameter one, however, $H_{15}$ is the smallest
that we know. A cover contains a unit segment, so one is the least
possible diameter. The sector has this diameter, and so do the earlier
covers cited above
\cite{gerriets-poole1974,norwood-poole-laidacker1992,norwood-poole2003,ploymaklam-wichiramala2018};
compare \cite[Figures~1 and~12]{movshovich2025}. Of these, only the
cover of \cite{ploymaklam-wichiramala2018} has smaller area than
$H_{15}$, and it is not convex. The other covers in
Table~\ref{tab:areas} that are smaller than the sector have diameter
greater than one: the hypotenuse of the scaled triangle has length
$1.0978\ldots$, and the other two were discussed above. We do not know whether a convex cover of smaller area than
$H_{15}$ has to be longer than one. The set conjectured in
\cite{movshovich2025} to have least area lies in the sector, so it has
diameter one, and by \cite{keller2026} its area would be at most
$0.2463322\ldots$. We return to this question in
Section~\ref{subsec:diameter}.

\begin{table}[htbp]
\centering
\caption{Areas of convex covers of unit arcs, rounded to seven decimal
places. The triangle is Wetzel's $30^\circ$--$60^\circ$--$90^\circ$
triangle scaled by the factor $1/1.0048$.}
\begin{tabular}{@{}llr@{}}
\toprule
Cover & Status & Area \\
\midrule
Wetzel's sector $\PiSector$ \cite{panraksa-wichiramala2021} & published & 0.2617994 \\
Scaled $30^\circ$--$60^\circ$--$90^\circ$ triangle \cite{wichiramala-panraksa2026} & preprint & 0.2609562 \\
$H_{15}$ & this paper & 0.2603092 \\
Triangle of Deng \cite{deng2026} & preprint & 0.2578836 \\
Quadrilateral of Keller \cite{keller2026} & preprint & 0.2463322 \\
\bottomrule
\end{tabular}

\label{tab:areas}
\end{table}

\subsection{The method}
\label{subsec:method}

Call a congruent copy of $\PiSector$ a \emph{unit sector}. Since only
one side of $\PiSector$ is trimmed, every unit sector contains two copies
of $H_{15}$, exchanged by the reflection in its bisector, and an arc in a
unit sector that has no congruent copy in $H_{15}$ lies in neither of
them. Such an arc meets both caps removed from that sector. For the trim
at $15^\circ$ these caps are disjoint (Figure~\ref{fig:trim}(b)), so the
arc has two distinct points in them, one on each side of the bisector.
The supporting wedges of \cite{panraksa-wichiramala2021} then reduce a
counterexample to a chain of contact points, points in the caps and
points at distance greater than one from an apex, in one of finitely many
orders that depend on two angles. It remains to show that every such
chain has length greater than one.

To each point $P_i$ of the chain we attach a vector $c_i$, with
$\sum_ic_i=0$. By Abel summation, $\sum_i\ip{c_i}{P_i}$ is at most the
length of the chain whenever every partial sum of the $c_i$, or
\emph{prefix}, lies in the unit disk (Lemma~\ref{lem:prefix}). We take
the $c_i$ at contact points normal to the supporting lines, so that their terms vanish, and the cap and
radial inequalities then show that the sum exceeds one. The same device
appears in Temerev and Doria's calibration for the golden gnomon
\cite[Lemma~4.2]{temerev-doria2026}, in Deng's balanced support
calibrations \cite[Lemma~2.7]{deng2026} and in Keller's proof
\cite[Section~4.3]{keller2026}. The first two bound suffix sums rather
than prefix sums, which comes to the same thing when the coefficients sum
to zero.

The configuration with an apex on the arc and Case~1 of
\cite{panraksa-wichiramala2021} are settled by hand. When the unit
sectors of both wedges contain the arc, elementary inequalities in one
angle suffice (Appendix~\ref{sec:capcap-circles}). Interval arithmetic is
needed only when the arc leaves the unit sector of one wedge through its
circular boundary and meets both caps of the other:
there we check prefix norms by rational interval arithmetic on finitely
many rectangles of the parameter domain.

What we take from \cite{panraksa-wichiramala2021} is stated as the Wedge
Lemma in Section~\ref{sec:reduction}. Apart from this and the prefix bound, the
proof has three new ingredients. First, a unit sector that contains the
arc contains two copies of $H_{15}$, which turns the single condition
that the arc has no copy in $H_{15}$ into two cap inequalities with
different normals. Second, a ray carrying two contacts has a single
normal, so its coefficient can be divided between the two contacts at
will; this makes one formula serve for all orientations in
Proposition~\ref{prop:case2-geometric}, and it supplies the parameters $a$
and $c$ of Section~\ref{subsec:balanced-normals}. Third,
Lemma~\ref{lem:boundary-prefixes} deals with two prefixes whose norm
reaches one, where interval arithmetic cannot decide.

The trim at $15^\circ$ was found by experiment. In a simplified model, a
chain consists of three points that alternate between the two radii of
one unit sector, three more that alternate starting on the other radius,
and one point of each closed cap of the trim $H_t$ of
Section~\ref{sec:trim}; there are
$8!/(3!\,3!)=1120$ orders that keep the order within each triple. For
each order we minimized the length of the chain with the local optimizer
L-BFGS-B, from one starting point.
For $t=16^\circ$ the shortest chain found was $O,\pt{14},\pt{16}$, of
length $1+2\sin1^\circ=1.0349\ldots$, and for $t=15.5^\circ$, $17^\circ$
and $20^\circ$ it was likewise $O,\pt{30-t},\pt t$, of length
$1+2\sin(t-15^\circ)$. For $t=15^\circ$ this chain is the radius
$[O,\pt{15}]$, of length one, whose end lies in both closed caps; it is
the equality case of Lemma~\ref{lem:apex}. For $t=12^\circ$ and
$14^\circ$ the optimizer found chains shorter than one, and this is why
we trim at $15^\circ$. The coefficients $\lambda,a,c$ of
Section~\ref{subsec:balanced-normals} were also found numerically. In
Case~2, for instance, a second-order cone program minimized the largest
prefix norm over nine sample points of a rectangle; the solution was
rounded to nine decimal places and kept if it passed the exact check
described in that section, and otherwise the rectangle was bisected. The
proof does not use these computations; it uses only the coefficients
they produced, which are checked in exact arithmetic. The two programs
are in the repository named at the end of the paper.

Section~\ref{sec:trim} sets up the trims $H_t$, Section~\ref{sec:local}
proves the general length bounds, Section~\ref{sec:reduction} reduces a
counterexample to finitely many ordered chains, and
Section~\ref{sec:proof} completes the proof.
Section~\ref{sec:threshold} is about deeper trims: it reduces the
inclusion-minimality of $H_{15}$ to two families of cuts, shows that
$H_t$ is not a cover for $t\le3.6^\circ$, and asks for the least area of
a convex cover of diameter one.

We write $R_\theta$ for the rotation through $\theta$, $\refl$ for the
reflection in the line through $O$ at angle $15^\circ$, and
$\ell(\gamma)$ for the length of an arc $\gamma$. For points $X,Y$ of a
simple arc $\gamma$ with parameters $\tau_X,\tau_Y$ we write $X\preceq Y$
if $\tau_X\le\tau_Y$, and $X\prec Y$ if $\tau_X<\tau_Y$.

\section{The trims and their caps}
\label{sec:trim}

\begin{definition}
\label{def:trim}
For $0^\circ<t<30^\circ$ put
\[
  m_t=\frac{t+30^\circ}{2},
  \qquad
  \delta_t=\frac{30^\circ-t}{2},
\]
and define the \emph{one-sided trim}
\begin{equation}
  H_t=\PiSector\cap
  \{x:\ip{x}{\pt{m_t}}\le\cos\delta_t\}.
  \label{eq:Ht}
\end{equation}
\end{definition}

The boundary line of the half-plane in \eqref{eq:Ht} passes through
$\pt t$ and $\pt{30}$, so $H_t$ keeps the circular boundary from $\pt0$
to $\pt t$ and replaces the arc from $\pt t$ to $\pt{30}$ by the chord
$[\pt t,\pt{30}]$. As the intersection of the compact convex sector with
a closed half-plane, $H_t$ is compact and convex. It is the union of the
circular sector $\Sigma_t$ and the triangle $\Delta_t$, where
\[
  \Sigma_t=\{r\pt\phi:0\le r\le1,\ 0^\circ\le\phi\le t\},
  \qquad
  \Delta_t=\conv\{O,\pt t,\pt{30}\},
\]
and these two pieces have disjoint interiors. To see that
$H_t=\Sigma_t\cup\Delta_t$, note that for $0^\circ\le\phi\le t$ the whole
radial segment at angle $\phi$ satisfies the half-plane condition in
\eqref{eq:Ht}, because
$\phi-m_t\in[-m_t,-\delta_t]$ gives
\begin{equation}
  \max_{0\le r\le1}\ \ip{r\pt\phi}{\pt{m_t}}
  =\cos(\phi-m_t)\le\cos\delta_t ,
  \label{eq:radial-in}
\end{equation}
while every point of $H_t$ with angle in $[t,30^\circ]$ lies in the closed
half-plane bounded by the chord line and containing $O$, hence in
$\Delta_t$.

\begin{proposition}[Area, monotonicity and nesting]
\label{prop:area}
For $0^\circ<t<30^\circ$,
\begin{equation}
  \area(H_t)=\frac{t\pi}{360}+\frac12\sin(30^\circ-t),
  \qquad
  \frac{d}{dt}\area(H_t)=\frac{\pi}{360}\bigl(1-\cos(30^\circ-t)\bigr)>0 .
  \label{eq:area-formula}
\end{equation}
The family is strictly nested: $H_s\subsetneq H_t$ for
$0^\circ<s<t<30^\circ$. In particular,
\[
  \area(H_{15})
  =
  \frac{\pi}{24}+\frac12\sin15^\circ
  =
  \frac{\pi}{24}+\frac{\sqrt6-\sqrt2}{8}
  =0.2603092\ldots .
\]
\end{proposition}

\begin{proof}
The decomposition $H_t=\Sigma_t\cup\Delta_t$ into pieces with disjoint
interiors gives the area formula, since $\area(\Sigma_t)=t\pi/360$ and
$\area(\Delta_t)=\frac12\sin(30^\circ-t)$. Differentiating with respect
to $t$, including the conversion from degrees to radians, gives the
displayed derivative, which is positive for $0^\circ<t<30^\circ$.

For nesting, let $s<t$. Clearly $\Sigma_s\subseteq\Sigma_t$. The triangle
$\Delta_s$ has vertices $O,\pt s,\pt{30}$; here $O,\pt{30}\in\Delta_t$ and
$\pt s\in\Sigma_t$, so convexity of $H_t$ gives $\Delta_s\subseteq H_t$
and hence $H_s\subseteq H_t$. The inclusion is strict because
$\pt t\in H_t$ while $\ip{\pt t}{\pt{m_s}}=\cos(t-m_s)>\cos\delta_s$, as
$|t-m_s|<\delta_s$ is equivalent to $s<t<30^\circ$.
\end{proof}

Every member of the family contains the triangle
\begin{equation}
 \conv\{O,\pt0,\pt{30}\}\subseteq H_t,
 \label{eq:triangular-core}
\end{equation}
since its three vertices satisfy the half-plane condition and $H_t$ is
convex.

\subsection{Reflection and caps}
\label{sec:chirality}

The reflection $\refl$ satisfies $\refl(\pt\phi)=\pt{30-\phi}$. The reflected trim $\refl(H_t)$
keeps the circular boundary for the angles in $[30^\circ-t,30^\circ]$ and
replaces the lower part of the arc by its chord.

\begin{proposition}[Mirror decomposition]
\label{prop:mirror-union}
For $15^\circ\le t<30^\circ$,
\[
  H_t\cup \refl(H_t)=\PiSector.
\]
For $t=15^\circ$ the two retained circular arcs meet exactly at
$\pt{15}$; for $t>15^\circ$ they overlap on the arc of angles
$[30^\circ-t,\,t]$. For $t<15^\circ$ the union is a proper subset of
$\PiSector$.
\end{proposition}

\begin{proof}
By \eqref{eq:radial-in}, $H_t$ contains the full radial segment at every
angle $\phi\in[0^\circ,t]$, and hence $\refl(H_t)$ contains the full radial
segment at every angle in $[30^\circ-t,30^\circ]$. If $t\ge15^\circ$ these
two intervals cover $[0^\circ,30^\circ]$, and every point of $\PiSector$
lies on some radial segment, so the union is $\PiSector$. The statements
about the arcs follow by intersecting the intervals $[0^\circ,t]$ and
$[30^\circ-t,30^\circ]$. If $t<15^\circ$, the points $\pt\phi$ with
$t<\phi<30^\circ-t$ belong to neither trim.
\end{proof}

\begin{figure}[htbp]
\centering
\includegraphics[width=0.92\linewidth]{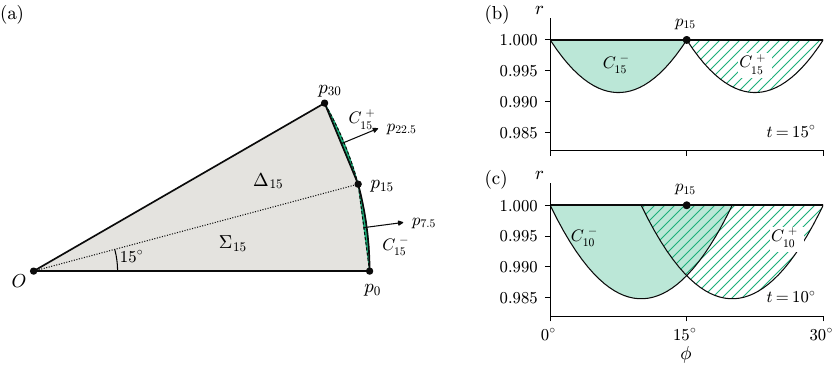}
\caption{(a) The trimmed sector $H_{15}$ (gray) is the union of the sector
$\Sigma_{15}$ and the triangle $\Delta_{15}=\conv\{O,\pt{15},\pt{30}\}$,
which meet along the dotted segment $[O,\pt{15}]$; the chord
$[\pt{15},\pt{30}]$ replaces the dashed arc of Wetzel's sector
$\PiSector$. The thin green regions are the cap $C^+_{15}$, cut off by
this chord, and the cap $C^-_{15}$, cut off by the chord
$[\pt0,\pt{15}]$ of $\refl(H_{15})$; their greatest width is
$1-\rho<0.009$. The arrows are the unit normals $\pt{22.5}$ and
$\pt{7.5}$ of the two chords.
(b) The open caps $C^-_{15}$ (shaded) and $C^+_{15}$ (hatched), removed
by $\refl(H_{15})$ and by $H_{15}$, drawn in the $(\phi,r)$-plane with the
$r$-axis stretched; their closures meet only at $\pt{15}$. (c) The caps
$C^\pm_{10}$ overlap for $10^\circ<\phi<20^\circ$, and both contain
$\pt{15}$.}
\label{fig:trim}
\end{figure}

The parts of $\PiSector$ removed by $\refl(H_t)$ and by $H_t$ are the
\emph{lower} and \emph{upper open caps}
\[
  C_t^-=
  \bigl\{x\in\PiSector:\ip{x}{\pt{\delta_t}}>\cos\delta_t\bigr\},
  \qquad
  C_t^+=
  \bigl\{x\in\PiSector:\ip{x}{\pt{m_t}}>\cos\delta_t\bigr\},
\]
respectively, which are relatively open in $\PiSector$; the \emph{closed caps}
$\overline C_t^\pm$ are defined with $\ge$ in place of $>$ and are
compact. A point $r\pt\phi$ of $C_t^+$ has $r\le1$, hence
$\cos(\phi-m_t)>\cos\delta_t$ and $t<\phi<30^\circ$; likewise the points
of $C_t^-$ have angles in $(0^\circ,30^\circ-t)$. For $t\ge15^\circ$ the
open caps are therefore disjoint. For $t=15^\circ$ the two caps are
exchanged by $\refl$, and every point of either cap has norm greater
than
\begin{equation}
 \rho=\cos7.5^\circ,
 \label{eq:rho}
\end{equation}
the distance from $O$ to each chord line. For $t<15^\circ$ the open caps
overlap, and both contain $\pt{15}$ (Figure~\ref{fig:trim}(c)).

For an apex $W$ and a direction $\theta$ write
\[
  \Pi(W,\theta)=W+R_\theta\PiSector
  =
  \{W+r\pt{\theta+\phi}:0\le r\le1,\ 0^\circ\le\phi\le30^\circ\},
\]
and abbreviate this unit sector to $\Pi_W$ when $\theta$ is clear. The two placements
of $H_{15}$ in $\Pi_W$ are
\[
H_W^+
=
\Pi_W\cap
\{X:\ip{X-W}{\pt{\theta+22.5}}\le \rho\},
\qquad
H_W^-
=
\Pi_W\cap
\{X:\ip{X-W}{\pt{\theta+7.5}}\le \rho\},
\]
the first removing the upper cap and the second the lower one. If
$\gamma\subseteq\refl(H_{15})$, then $\refl(\gamma)\subseteq H_{15}$; so an
arc that lies in either placement has a congruent copy in $H_{15}$.
Suppose that a unit arc $\gamma$ lies in $\Pi_W$ but in neither
placement. Since $H_W^\pm$ differ from $\Pi_W$ by one half-plane
condition each, there are points $X_W^+,X_W^-\in\gamma$ with
\begin{equation}
  \ip{X_W^+-W}{\pt{\theta+22.5}}>\rho,
  \label{eq:upper-cap}
\end{equation}
\begin{equation}
  \ip{X_W^--W}{\pt{\theta+7.5}}>\rho.
  \label{eq:lower-cap}
\end{equation}
We call such points \emph{witnesses}, and we always take a witness to be
a point of $\gamma$ at which the affine function on the left of
\eqref{eq:upper-cap} or \eqref{eq:lower-cap} attains its maximum over
$\gamma$.

Cap inequalities are unchanged by isometries applied to the witness, the
apex and the normal together: if $g(x)=Gx+a$ with $G$ orthogonal, then
\begin{equation}
 \ip{g(X)-g(W)}{Gq}=\ip{X-W}{q}
 \label{eq:transport}
\end{equation}
for every unit vector $q$.

\section{Length bounds from coefficients that sum to zero}
\label{sec:local}

\begin{lemma}[Prefix bound]
\label{lem:prefix}
Let $P_0,\ldots,P_N$ be points of an arc in parameter order, and let
$c_0,\ldots,c_N\in\R^2$ have sum zero. If every proper prefix
$S_k=\sum_{i=0}^k c_i$ $(0\le k<N)$ has norm at most $M$, where $M>0$,
then the length $L=\sum_{i=1}^N|P_i-P_{i-1}|$ of the polygonal chain
$P_0P_1\cdots P_N$ satisfies
\begin{equation}
 \Phi:=\sum_{i=0}^N\ip{c_i}{P_i}
 \le ML.
 \label{eq:prefix-radius}
\end{equation}
\end{lemma}

\begin{proof}
Since the $c_i$ sum to zero, Abel summation gives
\[
 \sum_{i=0}^N\ip{c_i}{P_i}
 =\sum_{i=0}^{N-1}\ip{S_i}{P_i-P_{i+1}},
\]
and each term on the right is at most $M|P_i-P_{i+1}|$.
\end{proof}

If the indices split into blocks whose coefficients each sum to zero,
then adding a fixed vector to every point of one block does not change
$\Phi$. We use this to move an apex to the origin.

\begin{remark}[Unfolding]
\label{rem:unfolding}
The reflection arguments of \cite{panraksa-wichiramala2021} are a special
case of Lemma~\ref{lem:prefix}. Choose for each segment $[P_{i-1},P_i]$ an
isometry $g_i(x)=G_ix+a_i$,
$1\le i\le N$, such that consecutive transformed segments meet,
$g_i(P_i)=g_{i+1}(P_i)$ for $1\le i<N$. For every unit vector $n$, putting
$v_i=G_i^{\mathsf T}n$ gives
\[
 \sum_{i=1}^N|P_i-P_{i-1}|
 \ge\sum_{i=1}^N\ip{v_i}{P_i-P_{i-1}}
 =\ip{n}{g_N(P_N)-g_1(P_0)},
\]
because each transformed segment keeps its length and the transformed
displacements telescope. The coefficients
\[
 c_0=-v_1,\qquad c_i=v_i-v_{i+1}\ (1\le i<N),\qquad c_N=v_N
\]
sum to zero, and their proper prefixes $\sum_{i=0}^k c_i=-v_{k+1}$
$(0\le k<N)$ have norm one. If the next piece is obtained by reflecting
the preceding one in the image of a supporting line
$\{x:\ip{x}{\nu_i}=h_i\}$ through $P_i$, with unit normal $\nu_i$, then
$v_{i+1}$ is the reflection of $v_i$ in the line through $0$ parallel to
that supporting line, so
\[
 c_i=2\ip{v_i}{\nu_i}\nu_i,\qquad
 \ip{c_i}{P_i}=2\ip{v_i}{\nu_i}h_i.
\]
This is why we take the coefficients at contact points normal to the
supporting lines. Lemma~\ref{lem:prefix} is more flexible than unfolding,
since it allows prefixes of norm less than one.
\end{remark}

\subsection{End wedges}
\label{subsec:endwedge}

\begin{figure}[htbp]
\centering
\includegraphics[width=0.92\linewidth]{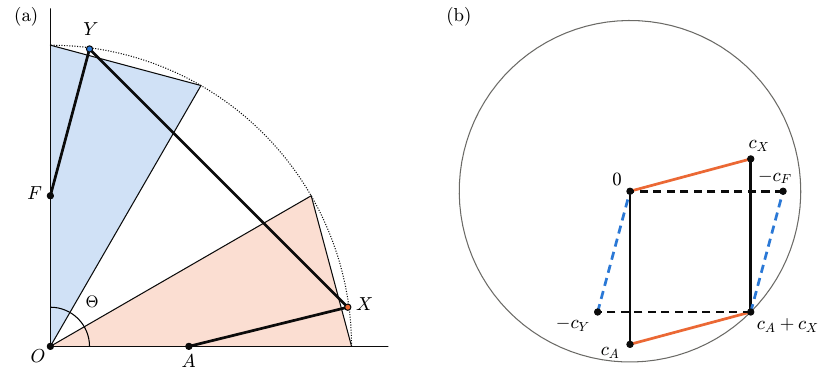}
\caption{Lemma~\ref{lem:end-wedge-tail-orders} for $\Theta=90^\circ$.
(a) The triangles $\conv\{O,\pt0,\pt{30}\}$ (orange) and
$\conv\{O,\pt{60},\pt{90}\}$ (blue), the points $A=0.46\,\pt0$, $X=0.996\,\pt{7.5}$,
$Y=0.996\,\pt{82.5}$ and $F=0.5\,\pt{90}$, and the chain $A,X,Y,F$; the
dotted arc is part of the unit circle; $X$ and $Y$ lie outside the shaded
triangles, as the lemma requires. (b) The vertices other than $0$ of the
parallelograms $0,c_A,c_A+c_X,c_X$ (solid) and $0,-c_Y,c_A+c_X,-c_F$
(dashed) are the proper prefixes of the four orders; both parallelograms
lie in the unit disk. Each side is a translate of one coefficient: $c_X$
in orange, $c_Y$ in blue, $c_A$ and $c_F$ in black.}
\label{fig:end-wedge}
\end{figure}

Let $60^\circ\le\Theta\le120^\circ$ and
$\Lambda_\Theta=\{r\pt\phi:r\ge0,\ 0^\circ\le\phi\le\Theta\}$. The triangle
\eqref{eq:triangular-core} lies in the unit disk and in the half-plane
$\ip{x}{\pt{7.5}}\le\rho$, as one checks at its three vertices. Hence a
point $X=r\pt\phi$ with $0^\circ\le\phi\le30^\circ$ and either $r>1$ or
$\ip{X}{\pt{7.5}}>\rho$ lies outside this triangle; since
$0^\circ\le\phi\le30^\circ$, this means
\begin{equation}
 \ip{X}{\pt{15}}>\cos15^\circ.
 \label{eq:end-chord-escape}
\end{equation}
Reflection in the bisector of $\Lambda_\Theta$ gives
$\ip{Y}{\pt{\Theta-15}}>\cos15^\circ$ for a point $Y=r\pt\phi$ with
$\Theta-30^\circ\le\phi\le\Theta$ and either $r>1$ or
$\ip{Y}{\pt{\Theta-7.5}}>\rho$. Only $\Theta=60^\circ$ is needed below
(Proposition~\ref{prop:case1-reflection}). We allow
$60^\circ\le\Theta\le120^\circ$ because the proof is the same and the
picture is clearer for a wider wedge (Figure~\ref{fig:end-wedge}).

\begin{lemma}[End wedges]
\label{lem:end-wedge-tail-orders}
Let $A$ and $F$ be points of the two rays that bound $\Lambda_\Theta$,
$A$ on the ray through $\pt0$ and $F$ on the ray through $\pt\Theta$. Let
$X=r\pt\phi$ with $0^\circ\le\phi\le30^\circ$ and either $r>1$ or
$\ip{X}{\pt{7.5}}>\rho$, and let $Y=r'\pt\psi$ with
$\Theta-30^\circ\le\psi\le\Theta$ and either $r'>1$ or
$\ip{Y}{\pt{\Theta-7.5}}>\rho$. Then each of the chains
\[
 A,X,Y,F;\qquad X,A,Y,F;\qquad
 A,X,F,Y;\qquad X,A,F,Y
\]
has length greater than $2\sin(\Theta/2)\ge1$.
\end{lemma}

\begin{proof}
Put
\[
 w=\frac{\sin(\Theta/2)}{\cos15^\circ},
 \qquad w'=\frac{\cos(\Theta/2-15^\circ)}{\cos15^\circ},
\]
and assign the coefficients
\[
 c_A=-w'\pt{90},\quad c_X=w\pt{15},\quad
 c_Y=w\pt{\Theta-15},\quad c_F=w'\pt{90+\Theta}.
\]
They sum to zero, and $c_A+c_X=(\sin(\Theta/2),-\cos(\Theta/2))$ is a
unit vector; for these identities use
$w'=\cos(\Theta/2)+\sin(\Theta/2)\tan15^\circ$ and reflection in the
bisector of $\Lambda_\Theta$. Since $30^\circ\le\Theta/2\le60^\circ$, we
have $|c_A|=|c_F|=w'\le1$ and $|c_X|=|c_Y|=w<1$. In each of the four
orders the first proper prefix is
$c_A$ or $c_X$, the second is $c_A+c_X$, and the third is $-c_Y$ or
$-c_F$, so all of them lie in the unit disk (Figure~\ref{fig:end-wedge}).
The support terms $\ip{c_A}{A}$ and $\ip{c_F}{F}$ vanish, so
Lemma~\ref{lem:prefix} and \eqref{eq:end-chord-escape} give, for the
length $L$ of the chain,
\[
 L\ge w\bigl(\ip{X}{\pt{15}}+\ip{Y}{\pt{\Theta-15}}\bigr)
   >2w\cos15^\circ=2\sin(\Theta/2). \qedhere
\]
\end{proof}

\subsection{The apex}
\label{subsec:apex}

Two cap normals symmetric about a unit vector $m$ add up to a multiple
of $m$, and this gives a single bound that holds whatever the order of
the apex and the two witnesses.

\begin{lemma}[The apex]
\label{lem:apex}
Let $m$ be a unit vector, let $q^-,q^+$ be its rotations through
$-\alpha$ and $\alpha$, where $0\le\alpha<60^\circ$, and put
$\kappa_\alpha=1/(2\cos\alpha)$. Let $h\in\R$, and let a rectifiable arc
$\gamma$ contain $O$, $X$ and $Y$. Put $s^-=\ip{X}{q^-}-h$ and
$s^+=\ip{Y}{q^+}-h$. Whatever the parameter order of $O,X,Y$,
\begin{equation}
 \ell(\gamma)\ge \frac{h}{\cos\alpha}+\kappa_\alpha(s^-+s^+).
 \label{eq:apex-identity}
\end{equation}
If $h=\cos\alpha$, $s^-=s^+=0$ and $\ell(\gamma)=1$, then $X=Y=m$ and
$\gamma$ is the segment $[O,m]$. In particular, an arc in $\PiSector$ that
contains $O$ and a point of each open cap $C_{15}^-,C_{15}^+$ has length
greater than one; for the closed caps the length is at least one, with
equality only for the radius $[O,\pt{15}]$.
\end{lemma}

\begin{proof}
Give $O$, $X$, $Y$ the coefficients $-m$, $\kappa_\alpha q^-$ and
$\kappa_\alpha q^+$. They sum to zero because $q^-+q^+=2\cos\alpha\, m$,
and their norms are $1,\kappa_\alpha,\kappa_\alpha$, where
$\kappa_\alpha<1$. In any order of three points a proper prefix is one
coefficient or minus another, so its norm is at most one, and
Lemma~\ref{lem:prefix} bounds $\ell(\gamma)$ below by
\[
 \kappa_\alpha\bigl(\ip{X}{q^-}+\ip{Y}{q^+}\bigr)
 =\frac{h}{\cos\alpha}+\kappa_\alpha(s^-+s^+).
\]
Suppose that $h=\cos\alpha$, $s^-=s^+=0$ and $\ell(\gamma)=1$, so that the
displayed quantity equals one. Then $O$ is not the middle point:
otherwise both proper prefixes would have norm $\kappa_\alpha$, and
Lemma~\ref{lem:prefix} with $M=\kappa_\alpha$ would give
$1\le\kappa_\alpha\ell(\gamma)=\kappa_\alpha<1$. If $O$ is an endpoint,
the edge joining $X$ and $Y$ is paired with a prefix of norm
$\kappa_\alpha<1$, so it has length zero; the other edge attains its
bound only along $m$, and the cap equations give $X=Y=m$. Equality of
chord length and arc length then gives $\gamma=[O,m]$. For the last
statement take $(\alpha,h,m)=(7.5^\circ,\rho,\pt{15})$: then $s^\pm$ are
positive for the open caps and nonnegative for the closed caps, and if
the length is one in the closed case, $s^-=s^+=0$.
\end{proof}

\section{Reduction to finitely many ordered chains}
\label{sec:reduction}

By \cite[Corollary~5]{wetzel-wichiramala2010}, a compact set in the plane
that contains a congruent copy of every simple polygonal unit arc
contains a congruent copy of every unit arc, simple or not; there, as
here, copies may be reflected. Since $H_{15}$ is compact, it suffices to
cover simple polygonal unit arcs.

For these we use the wedges of \cite{panraksa-wichiramala2021}. A
\emph{supporting wedge} of an arc $\gamma$ is a wedge of angle $30^\circ$
that contains $\gamma$ and whose two rays both meet $\gamma$. When two
supporting wedges with apices $U$ and $V$ are given, we choose
coordinates so that $U=O$, let $U_0$ and $U_1$ be the rays at $U$ of
directions $0^\circ$ and $30^\circ$, and let $V_0$ and $V_1$ be the rays
at $V$ of directions $\vartheta$ and $\vartheta+30^\circ$.

\begin{wedgelemma}[Panraksa and Wichiramala \cite{panraksa-wichiramala2021}]
Let $\gamma$ be a simple polygonal arc whose convex hull has nonempty
interior. Then either some supporting wedge of $\gamma$ has its apex on
$\gamma$, or $\gamma$ has two supporting wedges, with apices $U$ and $V$
not on $\gamma$, such that the following hold after reversing the
parametrization and reflecting the plane, if necessary.
\begin{enumerate}
\item[(a)] Each apex lies in the other wedge.
\item[(b)] The rays of each wedge carry three distinct points of
$\gamma$, with parameters $\tau_1<\tau_2<\tau_3$, the second on one ray
and the other two on the other ray. In Case~1 the triples are $A,B,D$ and
$A,C,D$; in Cases~2--4 they consist of six points $A,\ldots,F$. The
letters follow the parameter order, and Table~\ref{tab:incidences} gives
the ray that carries each point. Points of different triples may
coincide.
\item[(c)] In Case~1 the rays $U_0$ and $V_1$ lie on one line, so that
$\vartheta=150^\circ$, and the order on this line is $U,A,D,V$. In
Cases~2--4 the two points on one ray occur, outward from the apex, in the
orders $U,A,D$ and $V,F,C$ in Case~2, and $U,A,C$ and $V,F,D$ in Cases~3
and~4.
\item[(d)] If $\gamma$ has length one, then $\gamma$ lies in at least one
of the unit sectors $\Pi(U,0^\circ)$ and $\Pi(V,\vartheta)$.
\end{enumerate}
\end{wedgelemma}

\begin{table}[htbp]
\centering
\caption{The rays that carry the points of part~(b) of the Wedge Lemma in
Cases~1--4; the
first two columns list one triple and the last two the other.}
\begin{tabular}{@{}lcccc@{}}
\toprule
Case & $U_0$ & $U_1$ & $V_0$ & $V_1$\\
\midrule
1&$A,D$&$B$&$C$&$A,D$\\
2&$A,D$&$B$&$C,F$&$E$\\
3&$A,C$&$B$&$E$&$D,F$\\
4&$A,C$&$B$&$D,F$&$E$\\
\bottomrule
\end{tabular}

\label{tab:incidences}
\end{table}

The existence of the wedges, the three points on each, and statement~(a)
are given by \cite[Lemma~2.1]{panraksa-wichiramala2021} for the angle
$30^\circ$. The four cases are those of the proof of
\cite[Theorem~3.1]{panraksa-wichiramala2021}, and Figures~3 and~4 of that
paper show the positions of the points, in particular the orders in~(c).
That proof takes points $X,Y\in\gamma$ with $|X-V|>1$ and $|Y-U|>1$ and
shows in each case that $\gamma$ is longer than one, which gives~(d).

\begin{figure}[htbp]
\centering
\includegraphics[width=0.92\linewidth]{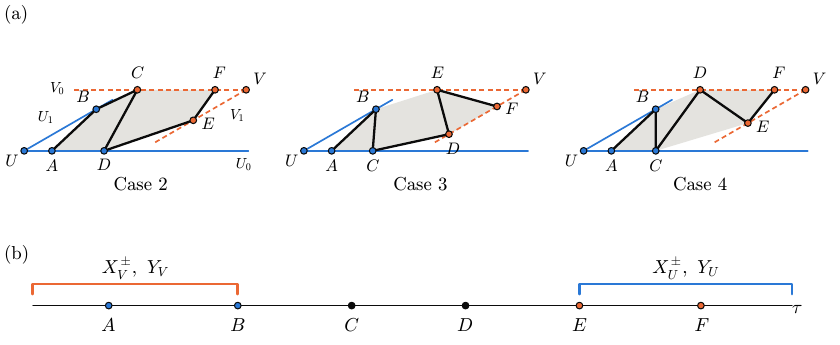}
\caption{(a) The contacts in Cases~2--4, drawn with $\vartheta=180^\circ$;
the rays are labeled in the first panel. The solid rays $U_0$ and $U_1$
(blue) issue from $U$, the dashed rays $V_0$ and $V_1$ (orange) from $V$,
and each contact has the color of the wedge whose ray carries it. The
polygonal arc through $A,\ldots,F$ is an example of an arc with these
contacts, and the gray polygon is its convex hull. (b) The parameter
order of the contacts along $\gamma$ and the ranges of
Lemma~\ref{lem:witness-bridge}(b): the $V$-witnesses and a point $Y_V$
farthest from $V$ occur no later than $B$ (orange bracket), and the
$U$-witnesses and a point $Y_U$ farthest from $U$ no earlier than $E$
(blue bracket). Here and in the later figures, blue refers to the wedge
at $U$ and orange to the wedge at $V$.}
\label{fig:canonical-orders}
\end{figure}

Now let $\gamma$ be a simple polygonal unit arc with no congruent copy in
$H_{15}$. Its convex hull has nonempty interior, because a simple arc of
length one on a line is a unit segment, and $H_{15}$ contains one. If a
supporting wedge has its apex $W$ on $\gamma$, then
$|X-W|\le\ell(\gamma)=1$ for every $X\in\gamma$, so $\gamma$ lies in the
unit sector of that wedge; since it lies in neither placement of
$H_{15}$, it contains points of both open caps, which
Lemma~\ref{lem:apex} excludes. So $\gamma$ has two supporting wedges as
in the Wedge Lemma; the reflection allowed there does no harm, because a
mirror image of $\gamma$ has no congruent copy in $H_{15}$ either. We
call the points of part~(b) the \emph{contacts}
(Figure~\ref{fig:canonical-orders}(a)). In Cases~2--4, writing
$\arg(V-U)=\xi$ and $\arg(U-V)=\vartheta+\zeta$ with
$\xi,\zeta\in[0^\circ,30^\circ]$, which is possible by (a), gives
$\vartheta=180^\circ+\xi-\zeta$, and hence
$150^\circ\le\vartheta\le210^\circ$. In Case~2 put $T=D$ and $P=C$; in
Cases~3 and~4 put $T=C$ and $P=D$. Then $A,B,T$ are the contacts on the
wedge at $U$, with $A,T\in U_0$ in the order $U,A,T$, and $P,E,F$ are
those on the wedge at $V$, with $F$ and $P$ on the same ray in the order
$V,F,P$.

For an apex $W\in\{U,V\}$ let $\theta_W$ be the direction of the first of
its two rays in counterclockwise order, so $\theta_U=0^\circ$ and
$\theta_V=\vartheta$. We say that $\gamma$
is in subcase $\sC_W$ (contained) if $\gamma\subseteq\Pi(W,\theta_W)$ and
in subcase $\sR_W$ (radial) otherwise. In subcase $\sC_W$, $\gamma$ has witnesses $X_W^+$ and
$X_W^-$ satisfying \eqref{eq:upper-cap} and \eqref{eq:lower-cap}. In
subcase $\sR_W$, a point $Y\in\gamma$ farthest from $W$ satisfies
$|Y-W|>1$, since $\gamma$ lies in the wedge at $W$ but not in $\Pi_W$.
Part~(d) of the Wedge Lemma excludes the subcase $\sR_U\sR_V$, so a
counterexample is in
one of the subcases
\begin{equation}
  \sC_U\sC_V,\qquad \sC_U\sR_V,\qquad \sR_U\sC_V,
  \label{eq:states}
\end{equation}
the last two being exchanged by a symmetry (see the proof of
Proposition~\ref{prop:mixed}). For a witness $X$ we use its cap slack
$s=\ip{X-W}{q}-\rho$, where $q$ is the normal in \eqref{eq:upper-cap} or
\eqref{eq:lower-cap}, and for a farthest point $Y$ its radial slack
$s=|Y-W|-1$. These slacks are positive for a counterexample.

\begin{lemma}[Ordered chains]
\label{lem:weak-order-enclosure}
Let the contacts of one case be listed in parameter order, and let
finitely many further labels of points of $\gamma$ be given. Then all
the labels admit an ordering $\omega=(\omega_1,\ldots,\omega_N)$ in which
the contacts keep their order and
$\tau_{\omega_1}\le\cdots\le\tau_{\omega_N}$. For the corresponding
points,
\[
  L_\omega
  :=
  \sum_{k=1}^{N-1}\bigl|P_{\omega_{k+1}}-P_{\omega_k}\bigr|
  \le \ell(\gamma).
\]
\end{lemma}

\begin{proof}
Order the classes of labels with equal parameters by their parameter
values; within a class, list the contacts in their given order and the
other labels arbitrarily. Labels with equal parameters name the same
point and add nothing to the sum. Each remaining term is at most the
length of the subarc between its two points, and the parameter intervals
of these subarcs overlap at most at endpoints.
\end{proof}

Deleting points from a chain in parameter order does not increase its
length, by the triangle inequality, so we may keep only some of the
witnesses. The next two results relate the position of a point on
$\partial\conv\gamma$ to its parameter.

\begin{lemma}[Crossing arcs]
\label{lem:alternating-boundary}
Let $K\subset\R^2$ be compact and convex with nonempty interior, and let
$Z_1,Z_2,Z_3,Z_4$ be distinct points of $\partial K$ such that the pairs
$\{Z_1,Z_2\}$ and $\{Z_3,Z_4\}$ alternate around $\partial K$. Then a
polygonal arc in $K$ from $Z_1$ to $Z_2$ and one from $Z_3$ to $Z_4$
intersect, even if their interiors meet $\partial K$.
\end{lemma}

\begin{proof}
Suppose that the two arcs are disjoint. They are compact, so they have
positive distance $d$. Insert a new vertex in the interior of each
segment, fix the four endpoints, and move each of the remaining
vertices a distance less than $d/3$ toward a fixed point of
$\operatorname{int}K$. The new arcs remain disjoint, since every point of
a new segment is within $d/3$ of the corresponding old segment, and by
convexity each new arc lies in $\operatorname{int}K$ except at its end
points. Removing polygonal loops gives two disjoint simple crosscuts of
$K$ with alternating endpoints, which is impossible by the Jordan curve
theorem; when $K$ is a polygon, the case used below, this is
\cite[Lemma~4.1.2(ii)]{diestel2025}.
\end{proof}

\begin{corollary}
\label{cor:localization}
Let $\gamma$ be a simple polygonal arc whose convex hull $K$ has nonempty
interior, and let
$A\prec B\prec T$ be distinct points of $\gamma\cap\partial K$. If $J$ is
the arc of $\partial K$ from $A$ to $B$ that does not contain $T$, then
every $X\in\gamma\cap J$ satisfies $X\preceq B$. The same holds with the
parameter order reversed.
\end{corollary}

\begin{figure}[htbp]
\centering
\includegraphics[width=0.45\linewidth]{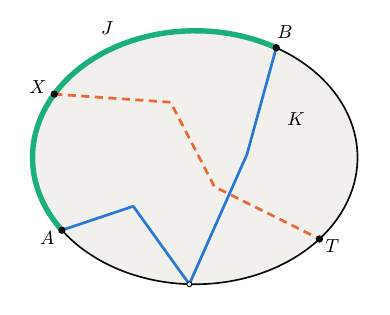}
\caption{Lemma~\ref{lem:alternating-boundary} and
Corollary~\ref{cor:localization}. The pairs $(A,B)$ and $(X,T)$ alternate
around $\partial K$, so a polygonal arc in $K$ from $A$ to $B$ (blue,
solid) and one from $X$ to $T$ (orange, dashed) must meet; here they do,
although the solid arc touches $\partial K$ at the open dot. Hence, for a
simple arc $\gamma$ with $A\prec B\prec T$, a point $X$ of $\gamma$ on the
arc $J$ (green, thick) of
$\partial K$ from $A$ to $B$ that does not contain $T$ satisfies
$X\preceq B$.}
\label{fig:localization}
\end{figure}

\begin{proof}
For $X=A$ or $X=B$ there is nothing to prove. Otherwise
$A,X,B,T$ are four distinct points in cyclic order on $\partial K$. If
$B\prec X$, the subarcs of $\gamma$ joining $A$ to $B$ and $T$ to $X$ have
disjoint parameter intervals, whatever the order of $T$ and $X$, so by
simplicity they are disjoint, contrary to
Lemma~\ref{lem:alternating-boundary} (Figure~\ref{fig:localization}).
\end{proof}

\begin{lemma}[Farthest points]
\label{lem:farthest}
Let $\gamma$ be an arc with convex hull $K$, and let $V\in\R^2$. Then
$\max_{z\in K}|z-V|=\max_{z\in\gamma}|z-V|$. If this maximum $r$ is
positive and is attained at $X\in\gamma$, then $X$ is the only point of $K$
at which $\ip{z}{n}$, with $n=(X-V)/r$, attains its maximum over $K$.
\end{lemma}

\begin{proof}
The function $z\mapsto|z-V|$ is convex, so its maximum over convex
combinations of points of $\gamma$ is at most its maximum over $\gamma$,
and $\gamma\subseteq K$ gives
the reverse inequality. For $z\in K$,
\[
 \ip{z-V}{n}\le |z-V|\le r=\ip{X-V}{n},
\]
and equality throughout forces $z-V=rn$, that is, $z=X$.
\end{proof}

\begin{lemma}[Localization]
\label{lem:witness-bridge}
Let $K=\conv\gamma$.
\begin{enumerate}
\item[(a)] In Case~1, let $X$ be a point of $\gamma$ maximizing
$\ip{x-V}{\pt{172.5}}$ if $\gamma\subseteq\Pi_V$, and a point of
$\gamma$ farthest from $V$ otherwise; let $Y$ be a point maximizing
$\ip{y}{\pt{7.5}}$ if $\gamma\subseteq\Pi_U$, and a point farthest from
$U$ otherwise. Then $X\preceq B$ and $C\preceq Y$.
\item[(b)] In Cases~2--4, both $V$-witnesses $X_V^\pm$ and a point of
$\gamma$ farthest from $V$ occur no later than $B$, and both $U$-witnesses
$X_U^\pm$ and a point farthest from $U$ occur no earlier than $E$
(Figure~\ref{fig:canonical-orders}(b)).
\item[(c)] Consequently, in subcase $\sR_U\sC_V$ of Cases~2--4 the points
$A,X_V^-,X_V^+,B,\ldots,F$ and a point $Y=r\pt\phi$ of $\gamma$
farthest from $U$, where $r>1$ and $0^\circ\le\phi\le30^\circ$, occur in
one of the twelve orders
\begin{equation}
 \{A,X_V^-,X_V^+\},B,C,D,E,Y,F
 \quad\text{or}\quad
 \{A,X_V^-,X_V^+\},B,C,D,E,F,Y,
 \label{eq:mixed-twelve}
\end{equation}
where $\{A,X_V^-,X_V^+\}$ stands for these three labels in any of their
six orders; the proper prefixes of these orders form $13$ distinct sets
of labels. In subcase $\sC_U\sC_V$, a $V$-witness
$X_V$ and a $U$-witness $X_U$ occur in one of the four orders
\begin{equation}
 (A,X_V\text{ or }X_V,A),B,C,D,E,
 (X_U,F\text{ or }F,X_U),
 \label{eq:capcap-four}
\end{equation}
whose proper prefixes form nine distinct sets of labels.
\end{enumerate}
\end{lemma}

\begin{proof}
(b) For the wedge at $U$, $A\prec B\prec T$ are the distinct contacts of
one triple, and the order on $U_0$ is $U,A,T$. The lines that contain
$U_1$ and $U_0$ support $K$, with outer normals $\pt{120}$ and
$\pt{270}$: the first at $B$, the second along a segment that contains
$A$ and $T$. As $\partial K$ is traversed counterclockwise from $B$, the
outer normal turns from $\pt{120}$ to $\pt{270}$, and the segment on the
line of $U_0$ is traversed away from $U$, so that $A$ is reached before
$T$. The part of
$\partial K$ from $B$ to $A$ is therefore the arc $J_U$ from $A$ to $B$
that does not contain $T$, and it contains every face of $K$ whose
normal has angle strictly between $120^\circ$ and $270^\circ$. The cap
normals at $V$ and the directions from $V$ to points of its wedge have
angles in $[\vartheta,\vartheta+30^\circ]\subseteq[150^\circ,240^\circ]$. A
$V$-witness maximizes such a functional over $\gamma$, hence over $K$, so
it lies on the face of $K$ exposed by that normal, and this face lies on
$J_U$; by Lemma~\ref{lem:farthest}, so does a point farthest from $V$.
Such a point $X$ satisfies $X\preceq B$: trivially if it is $A$ or $B$,
and by Corollary~\ref{cor:localization} otherwise.

For the wedge at $V$, $P\prec E\prec F$ are the distinct contacts of the other
triple, and the order on the ray that carries $F$ and $P$ is $V,F,P$. The
lines that contain $V_1$ and $V_0$ support $K$ with outer normals
$\pt{\vartheta+120}$ and $\pt{\vartheta+270}$. Going around $\partial K$
from $E$, counterclockwise in Cases~2 and~4 and clockwise in Case~3, one
finds in the same way that the arc of $\partial K$ from $E$ to $F$ that
does not contain $P$ contains every face whose normal has angle strictly
between $\vartheta+120^\circ$ and $\vartheta+270^\circ$. The directions
from $U$ to points of its wedge, lifted from $[0^\circ,30^\circ]$ to
$[360^\circ,390^\circ]$, lie strictly between these bounds for
$150^\circ\le\vartheta\le210^\circ$, and so do the two cap normals at
$U$. The argument above, with the parameter order reversed, gives
$E\preceq Y$ for each $U$-witness $Y$ and for a point $Y$ farthest from
$U$. Only the
strict order of the three contacts within each triple was used, so the
argument applies also when the face exposed by a normal is an edge, or
when a witness coincides with a contact.

(a) The same argument applies in Case~1. The normal that defines $X$ has
angle $172.5^\circ$ or, if $X$ is a point farthest from $V$, an angle in
$[150^\circ,180^\circ]$; both lie strictly between $120^\circ$ and
$270^\circ$, so $X$ lies on the arc of $\partial K$ from $A$
to $B$ that does not contain $D$. At the other end the normal of angle
$7.5^\circ$, or a direction in $[0^\circ,30^\circ]$, lifted by $360^\circ$
lies strictly between $270^\circ$ and $420^\circ$, the range of normals
facing $V$, so
$Y$ lies on the arc from $C$ to $D$ that does not contain $A$.
Corollary~\ref{cor:localization}, applied to the triples
$A\prec B\prec D$ and $A\prec C\prec D$, the second with the parameter
order reversed, gives $X\preceq B$ and $C\preceq Y$, also when $B=C$ or
when $X$ or $Y$ is a contact.

(c) By (b), the labels $A,X_V^-,X_V^+$ precede $B$ in some order, and $Y$
and $F$ follow $E$ in some order; with Lemma~\ref{lem:weak-order-enclosure}
this gives \eqref{eq:mixed-twelve}. With one witness at each end the same
argument gives \eqref{eq:capcap-four}.
\end{proof}

From now on we identify a proper prefix with its set of labels and
write it in braces, as in $\{A,B,X_V\}$.

\begin{proposition}[Reduction]
\label{prop:certificate-criterion}
Suppose that, in each of Cases~1--4 and each subcase in \eqref{eq:states},
every simple polygonal counterexample $\gamma$ yields an ordered chain
$P_1,\ldots,P_n$ of its points, as in
Lemma~\ref{lem:weak-order-enclosure}, of length greater than one. Then
$H_{15}$ covers every unit arc.
\end{proposition}

\begin{proof}
By Lemma~\ref{lem:weak-order-enclosure} the chain has length at most
$\ell(\gamma)=1$, so no simple polygonal counterexample exists in these
cases; the apex case and the subcase $\sR_U\sR_V$ were excluded above.
Compactness and \cite[Corollary~5]{wetzel-wichiramala2010} give the claim
for all unit arcs.
\end{proof}

\section{Proof of Theorem~\ref{thm:H15}}
\label{sec:proof}

By Proposition~\ref{prop:certificate-criterion} it remains to bound the
length of the ordered chains of Lemma~\ref{lem:witness-bridge} from
below. Throughout this section
\[
 e_0=\pt{90},\qquad e_1=\pt{-60},\qquad \kappa=\frac1{2\rho};
\]
here $e_0$ and $e_1$ are unit normals to the rays $U_0$ and $U_1$, and
$\kappa=\kappa_{7.5^\circ}$ in the notation of Lemma~\ref{lem:apex}.

\subsection{Case 1}
\label{subsec:case1-reflection}

Here the two wedges share a supporting line. Reflecting the two ends of
the chain in the other two supporting lines makes the apices coincide,
and Lemma~\ref{lem:end-wedge-tail-orders} applies.

\begin{proposition}[Case 1]
\label{prop:case1-reflection}
Every simple polygonal unit arc in Case~1 of the Wedge Lemma lies in at
least one of the four placements $H_U^\pm$, $H_V^\pm$.
\end{proposition}

\begin{figure}[htbp]
\centering
\includegraphics[width=0.92\linewidth]{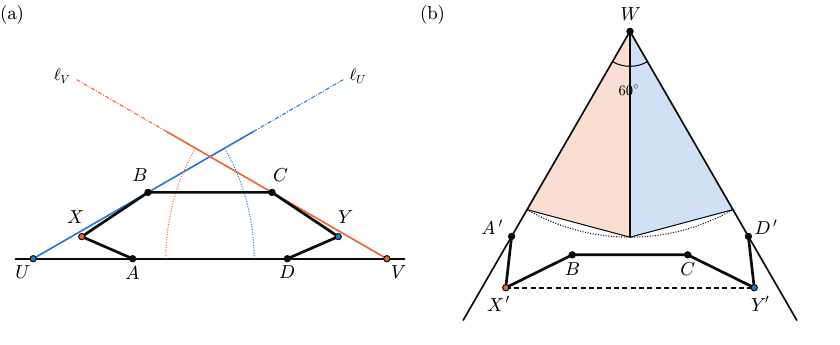}
\caption{Proposition~\ref{prop:case1-reflection}, drawn with $d=1.6$.
(a) The common supporting line through $U$ and $V$; the rays $U_1$ (blue)
and $V_0$ (orange), which make angles of $30^\circ$ with it, and their
lines $\ell_U$ and $\ell_V$ (continued dash-dotted); the unit circles
about $U$ and $V$ (dotted, in the same colors); the contacts $A,B,C,D$;
and the points $X,Y$ of Lemma~\ref{lem:witness-bridge}(a), here the
points of the chain farthest from $V$ and from $U$. (b) After the part of
the chain up to $B$ is reflected in $\ell_U$ and the part from $C$ in
$\ell_V$, the apices coincide at $W$ and the chain $A',X',B,C,Y',D'$ lies
in the $60^\circ$ wedge at $W$. The shaded end triangles and the unit
circle about $W$ (dotted) are those of
Lemma~\ref{lem:end-wedge-tail-orders} for $\Theta=60^\circ$, colored as
in Figure~\ref{fig:end-wedge}; the dashed segment shortens the chain to
$A',X',Y',D'$.}
\label{fig:case1-unfolding}
\end{figure}

\begin{proof}
Suppose that such an arc $\gamma$ lies in none of them. Choose coordinates so that
$U=(0,0)$ and $V=(d,0)$, with $d>0$, lie on the common supporting line;
the rays at $U$ have
directions $0^\circ$ and $30^\circ$, and those at $V$ have directions
$150^\circ$ and $180^\circ$. The contacts satisfy
$A\preceq B\preceq C\preceq D$, with $A,D$ on the common line, $B$ on the
ray $U_1$ and $C$ on the ray $V_0$ of direction $150^\circ$
(Figure~\ref{fig:case1-unfolding}(a)). Let $X$ and $Y$ be the points of
Lemma~\ref{lem:witness-bridge}(a). If $\gamma\subseteq\Pi_V$, then
$\ip{X-V}{\pt{172.5}}>\rho$, because $\gamma\not\subseteq H_V^+$;
otherwise $|X-V|>1$. Likewise $\ip{Y}{\pt{7.5}}>\rho$ or $|Y-U|>1$. By
that lemma, $X\preceq B$ and $C\preceq Y$.

Let $\ell_U$ and $\ell_V$ be the lines that contain $U_1$ and $V_0$. The
reflection of $V$ in $\ell_U$ and the reflection of $U$ in $\ell_V$ are
the same point,
\[
 W=\Bigl(\frac d2,\frac{\sqrt3\,d}{2}\Bigr).
\]
Reflect the part of the chain up to $B$ in $\ell_U$ and the part from $C$
in $\ell_V$. These reflections fix $B$ and $C$, so the reflected pieces
still form a chain of the same length. Let $A'$ and $X'$ be the images of
$A$ and $X$ under the first reflection, and $Y'$ and $D'$ those of $Y$
and $D$ under the second. Seen from $W$, the points $A'$ and $D'$ lie on
the rays of directions $240^\circ$ and $300^\circ$, and the directions of
$X'-W$ and $Y'-W$ lie in $[240^\circ,270^\circ]$ and
$[270^\circ,300^\circ]$, respectively. The reflections preserve distances
to the corresponding apex, and by \eqref{eq:transport} they take the cap
normals of directions $172.5^\circ$ and $7.5^\circ$ to the directions
$247.5^\circ$ and $292.5^\circ$. These are the normals $\pt{7.5}$ and
$\pt{\Theta-7.5}$ of
Lemma~\ref{lem:end-wedge-tail-orders} for $\Theta=60^\circ$, rotated
through $240^\circ$. Deleting $B$ and $C$ from the chain leaves one of the
chains
\[
 A',X',Y',D';\qquad X',A',Y',D';\qquad
 A',X',D',Y';\qquad X',A',D',Y',
\]
and Lemma~\ref{lem:end-wedge-tail-orders}, translated to $W$ and rotated
through $240^\circ$, shows that each of them is longer than
$2\sin30^\circ=1$ (Figure~\ref{fig:case1-unfolding}(b)). But the chain has
length at most $\ell(\gamma)=1$, a contradiction.
\end{proof}

\subsection{Both sectors contain the arc}
\label{subsec:capcap}

In subcase $\sC_U\sC_V$ of Cases~2--4 we keep one witness at each apex,
so that Lemma~\ref{lem:witness-bridge}(c) leaves the four orders
\eqref{eq:capcap-four}. In Case~2 a single formula for the coefficients
serves for all $\vartheta\in[150^\circ,210^\circ]$. This is possible
because the coefficient of a ray that carries two contacts can be
divided between them in any proportion without changing $\Phi$.

\begin{proposition}[Case 2]
\label{prop:case2-geometric}
In subcase $\sC_U\sC_V$ of Case~2, let $X_U=X_U^+$ and $X_V=X_V^+$ be the
upper-cap witnesses, and put
\[
 s_U=\ip{X_U}{\pt{22.5}}-\rho,\qquad
 s_V=\ip{X_V-V}{\pt{\vartheta+22.5}}-\rho.
\]
For each of the four orders in \eqref{eq:capcap-four}, the nonstrict
inequalities $s_U,s_V\ge0$ imply that the chain has length
\begin{equation}
 L\ge\frac{50}{49}\bigl(1+\kappa(s_U+s_V)\bigr)
 \ge\frac{50}{49}>1.
 \label{eq:case2-global-length}
\end{equation}
\end{proposition}

\begin{proof}
By Lemma~\ref{lem:witness-bridge}(c), $X_U$ and $X_V$ occur in one of the
orders \eqref{eq:capcap-four}. Put $b=2\cos15^\circ-1$, so that
$Q=e_0+b e_1=\kappa\pt{22.5}$, and assign the coefficients
\[
\begin{array}{c|c@{\qquad}c|c}
A&a e_0&C&cR_\vartheta e_0\\
B&-b e_1&E&-bR_\vartheta e_1\\
D&(-1-a)e_0&F&(-1-c)R_\vartheta e_0\\
X_U&Q&X_V&R_\vartheta Q
\end{array}
\qquad
 a=-\frac{2+\sin\vartheta}{3},\quad
 c=-\frac{1+\sin\vartheta}{3}.
\]
The coefficients at each apex sum to zero, so we may replace each point
$Z$ of the block at $V$ by $Z-V$ (see the paragraph after
Lemma~\ref{lem:prefix}). Every support term then vanishes, because each
support coefficient is normal to the ray through $U$ or $V$ that carries
its contact, and $\Phi=1+\kappa(s_U+s_V)$. By
Appendix~\ref{subsec:case2-global-bounds}, every proper prefix has norm
less than $49/50$ for $150^\circ\le\vartheta\le210^\circ$, and
\eqref{eq:prefix-radius} with $M=49/50$ gives
\eqref{eq:case2-global-length}.
\end{proof}

\begin{figure}[htbp]
\centering
\includegraphics[width=0.62\linewidth]{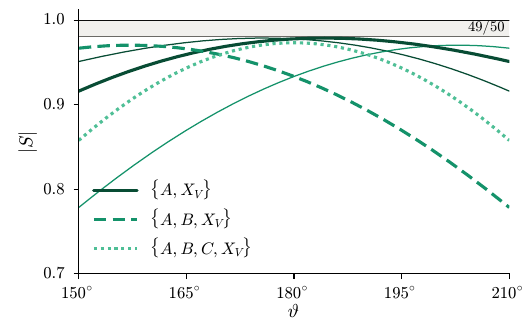}
\caption{The norms of the prefixes $\{A,X_V\}$ (solid), $\{A,B,X_V\}$
(dashed) and $\{A,B,C,X_V\}$ (dotted) of the coefficients in
Proposition~\ref{prop:case2-geometric}, for
$150^\circ\le\vartheta\le210^\circ$; the gray band lies between the
bound $49/50$ and $1$. The thin curves are the norms of
$\{A,\ldots,D,X_V\}$ and $\{A,\ldots,E,X_V\}$, which are those of the
second and the first prefix reflected in $\vartheta=180^\circ$; the
remaining four prefixes have norm at most $5/6$.}
\label{fig:case2-prefix-loci}
\end{figure}

By the symmetry $\vartheta\mapsto360^\circ-\vartheta$ of
Appendix~\ref{subsec:case2-global-bounds}, it suffices to bound three
prefix norms; Figure~\ref{fig:case2-prefix-loci} shows them.

\begin{proposition}[Cases 3 and 4]
\label{prop:case34}
In subcase $\sC_U\sC_V$ of Cases~3 and~4, let
\[
 (\phi_U,\phi_V)=
 \begin{cases}
 (22.5^\circ,7.5^\circ)&\text{if }150^\circ\le\vartheta\le180^\circ,\\
 (7.5^\circ,22.5^\circ)&\text{if }180^\circ\le\vartheta\le210^\circ,
 \end{cases}
\]
where either pair may be used at $\vartheta=180^\circ$. Let $X_U$ and
$X_V$ be the witnesses
for the normals $\pt{\phi_U}$ and $\pt{\vartheta+\phi_V}$, and let
$s_U,s_V$ be their slacks. Then, in each of the four orders
\eqref{eq:capcap-four}, the chain has length at least
$1+\kappa(s_U+s_V)>1$.
\end{proposition}

\begin{proof}
For the coefficients given in Appendix~\ref{subsec:case34-circles}, those
at each apex sum to zero and $\Phi=1+\kappa(s_U+s_V)$. On each of the
intervals $[150^\circ,180^\circ]$ and $[180^\circ,210^\circ]$ the
coefficients at $U$ are fixed and those at $V$ rotate with $\vartheta$, so
every prefix has the form $u+R_\vartheta v$ with $u,v$ fixed. For
$\vartheta$ in an interval $[\vartheta_1,\vartheta_2]$ with
$0<\vartheta_2-\vartheta_1<180^\circ$, the point $u+R_\vartheta v$ runs
along a circle of radius $|v|$ about $u$. Since
$|u+R_\vartheta v|^2=|u|^2+|v|^2+2\ip{u}{R_\vartheta v}$, its norm is
largest at an endpoint unless $R_\vartheta v$ points in the direction of
$u$ at an interior parameter, and then the largest norm is $|u|+|v|$; for
$u,v\ne0$ this happens exactly when $\det(R_{\vartheta_1}v,u)>0$ and
$\det(u,R_{\vartheta_2}v)>0$. Of the nine proper prefixes, four have
constant norm at most $3/4$, and four move on such circles;
Appendix~\ref{subsec:case34-circles} shows that each of these four has
norm less than $999/1000$. The remaining prefix, $\{A,B,C,X_V\}$, is
$\kappa\pt{\vartheta+\phi_V}-\kappa\pt{\phi_U}$. The pairs
$(\phi_U,\phi_V)$ are exchanged at $\vartheta=180^\circ$ so that the angle
between $\pt{\phi_U}$ and $\pt{\vartheta+\phi_V}$ never exceeds
$165^\circ$, and then this prefix has norm at most
$2\kappa\sin82.5^\circ=1$ (Appendix~\ref{subsec:case34-circles}). Hence
every proper prefix lies in the unit disk, and Lemma~\ref{lem:prefix}
gives the bound.
\end{proof}

\subsection{One sector contains the arc}
\label{subsec:balanced-normals}

\begin{figure}[!tb]
\centering
\includegraphics[width=0.92\linewidth]{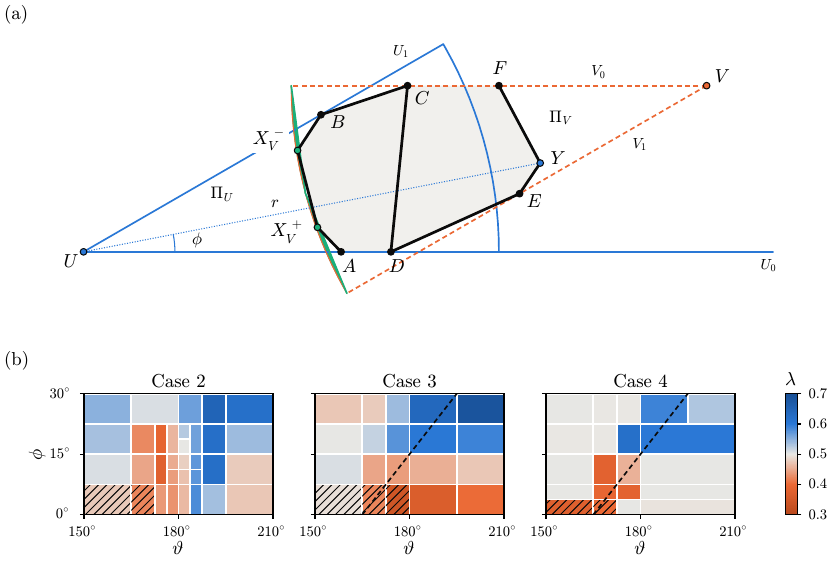}
\caption{Subcase $\sR_U\sC_V$. (a) A chain with the contacts of Case~2,
drawn with $\vartheta=180^\circ$ and $\phi=11^\circ$: the unit sectors
$\Pi_U$ (blue) and $\Pi_V$ (orange, with dashed rays), the two caps of
$\Pi_V$ (green, between its circular arc and the two chords), the
witnesses $X_V^-$ and $X_V^+$ in these caps, and a point $Y=r\pt\phi$
farthest from $U$, which lies outside $\Pi_U$. The chain shown is much
longer than one; Proposition~\ref{prop:mixed} shows that every such chain
is longer than one. (b) The rectangles that cover the domain of
$(\vartheta,\phi)$ in Cases~2, 3 and~4, colored by the weight $\lambda$:
blue where the radial inequality at $U$ has the larger weight, orange
where the two cap inequalities at $V$ have. On the hatched rectangles
$a=-\sqrt3\lambda$, and Lemma~\ref{lem:boundary-prefixes}(i) is used. In
Cases~3 and~4 the prefix of Lemma~\ref{lem:boundary-prefixes}(ii) has
norm one on the dashed line $\vartheta=165^\circ+\phi$.}
\label{fig:mixed}
\end{figure}

In subcase $\sR_U\sC_V$ the chain has a point farthest from $U$ at one end
and both $V$-witnesses at the other (Figure~\ref{fig:mixed}(a)). We combine the radial and the cap
inequalities with a weight $\lambda$; as in Section~\ref{subsec:capcap},
the parameters $a$ and $c$ divide the coefficient of each ray that
carries two contacts between those two contacts.

\begin{lemma}[Coefficients for subcase $\sR_U\sC_V$]
\label{lem:common-mixed}
In subcase $\sR_U\sC_V$ of Cases~2--4, let $Y=r\pt\phi$ be a point of
$\gamma$ farthest from $U$, with $r>1$ and $0^\circ\le\phi\le30^\circ$, and
let $X_V^\pm$ be the two $V$-witnesses. Put
\[
 d_i=R_\vartheta e_i,
 \qquad q^-=\pt{7.5},\quad q^+=\pt{22.5},
 \qquad s_V^\pm=\ip{X_V^\pm-V}{R_\vartheta q^\pm}-\rho,
\]
let $\varepsilon=0$ in Cases~2 and~4 and $\varepsilon=1$ in Case~3, and
for $0<\lambda<1$ and $a,c\in\R$ set
\[
 \mu_0=-\lambda(\sin\phi+\sqrt3\cos\phi),\quad
 \mu_1=-2\lambda\cos\phi,\quad
 g=-2(1-\lambda)\cos15^\circ.
\]
Assign the coefficients
\[
\begin{array}{c|c@{\qquad}c|c}
 A&a e_0&P&c d_\varepsilon\\
 B&\mu_1 e_1&E&g d_{1-\varepsilon}\\
 T&(\mu_0-a)e_0&F&(g-c)d_\varepsilon\\
 Y&\lambda\pt\phi&X_V^\pm&(1-\lambda)\kappa R_\vartheta q^\pm .
\end{array}
\]
If all proper prefixes in one of the orders \eqref{eq:mixed-twelve} have
norm at most one, then the chain has length
\begin{equation}
 L\ge1+\lambda(r-1)+(1-\lambda)\kappa(s_V^-+s_V^+)>1.
 \label{eq:case2-mixed-functional}
\end{equation}
\end{lemma}

\begin{proof}
We have
\[
 \mu_0 e_0+\mu_1 e_1+\lambda\pt\phi=0,
 \qquad q^-+q^+=2\rho\pt{15},
 \qquad g(e_0+e_1)+(1-\lambda)\pt{15}=0.
\]
Hence the coefficients at each apex sum to zero. At $U$ the three support
coefficients are perpendicular to their rays, so their terms vanish. At
$V$ we replace each point $Z$ of the block by $Z-V$, which does not change
$\Phi$, and the support terms then vanish as well. What remains is
\[
 \lambda r+(1-\lambda)\kappa\bigl(\ip{X_V^--V}{R_\vartheta q^-}
           +\ip{X_V^+-V}{R_\vartheta q^+}\bigr)
 =1+\lambda(r-1)+(1-\lambda)\kappa(s_V^-+s_V^+),
\]
and Lemma~\ref{lem:prefix} gives the bound, which is strict because
$r>1$ and $s_V^\pm>0$.
\end{proof}

Two prefixes can reach norm one; we bound them separately.

\begin{lemma}[Two prefixes of norm at most one]
\label{lem:boundary-prefixes}
In the setting of Lemma~\ref{lem:common-mixed}:
\begin{enumerate}
\item[(i)] if $a=-\sqrt3\lambda$, the prefix $\{A,B,X_V^-,X_V^+\}$ is
$\lambda w+(1-\lambda)\pt{\vartheta+15}$, where
\[
 w=(-\cos\phi,\sqrt3(\cos\phi-1))
   =(2\cos\phi-1)(-\pt0)+2(1-\cos\phi)(-\pt{60}),
\]
and this prefix has norm at most one;
\item[(ii)] in Cases~3 and~4 the prefix $\{A,B,C,X_V^-,X_V^+\}$ is
$-\lambda\pt\phi+(1-\lambda)\pt{\vartheta+15}$; it has norm at most one,
with equality when $\vartheta=165^\circ+\phi$.
\end{enumerate}
\end{lemma}

\begin{proof}
(i) The part from $U$ is $a e_0+\mu_1 e_1=\lambda w$, and the part from
$V$ is $(1-\lambda)\kappa R_\vartheta(q^-+q^+)=(1-\lambda)\pt{\vartheta+15}$. For
$0^\circ\le\phi\le60^\circ$ the coefficients $2\cos\phi-1$ and
$2(1-\cos\phi)$ are nonnegative and sum to one, so $w$ is a convex
combination of the unit vectors $-\pt0$ and $-\pt{60}$ and lies in the
unit disk, and the prefix is a convex combination of $w$ and a unit
vector. (ii) In Cases~3 and~4, $T=C$, so the prefix contains the whole
block at $U$ except $Y$, and this part sums to $-\lambda\pt\phi$. The
prefix is therefore a convex combination of the unit vectors $-\pt\phi$
and $\pt{\vartheta+15}$, which coincide exactly when
$\vartheta+15^\circ=\phi+180^\circ$.
\end{proof}

For each of Cases~2--4 we cover the domain
$[150^\circ,210^\circ]\times[0^\circ,30^\circ]$ of $(\vartheta,\phi)$ by
finitely many closed rectangles, each with its own coefficients
$\lambda,a,c$ (Figure~\ref{fig:mixed}(b)); these are rational, except
that $a=-\sqrt3\lambda$ on some rectangles. A program checks that the rectangles cover the domain, that
$0<\lambda<1$, and that every proper prefix of each of the twelve orders
\eqref{eq:mixed-twelve} has squared norm at most one on its rectangle,
except the prefixes of Lemma~\ref{lem:boundary-prefixes}, which are
accepted only where the hypotheses of that lemma hold. To check
coverage, the program cuts the domain into cells along the endpoints of
all rectangles and checks that each cell lies in some rectangle. The norms are bounded by interval arithmetic
with rational endpoints: sine and cosine are enclosed by alternating
Taylor bounds with explicit remainders after quadrant reduction, $\pi$ by
Machin's formula, and the endpoints of the intervals are rounded outward
to $28$ decimal places, which keeps the denominators small. A prefix is
accepted on a rectangle when the upper bound for its squared norm is at
most one; otherwise the program bisects the longer side and checks both
halves, and a prefix still undecided after $28$ bisections in Case~2 or
$30$ in Cases~3 and~4 makes the check fail. The largest bisection depths
reached are $13$, $9$ and $10$ in Cases~2, 3 and~4, respectively. In
Cases~3 and~4, Lemma~\ref{lem:boundary-prefixes}(ii) is used on every
rectangle, and Lemma~\ref{lem:boundary-prefixes}(i) settles the prefix
$\{A,B,X_V^-,X_V^+\}$ on two rectangles in Case~2, three in Case~3 and two
in Case~4. Table~\ref{tab:mixed-counts} gives the size of the
computation. The whole check, including the programs for
Appendix~\ref{sec:capcap-circles}, takes about a minute or less on a
laptop.

\begin{table}[htbp]
\centering
\caption{The computation for subcase $\sR_U\sC_V$: the number of rectangles
covering $[150^\circ,210^\circ]\times[0^\circ,30^\circ]$, and the number
of pairs (rectangle, prefix) settled by interval bounds and by
Lemma~\ref{lem:boundary-prefixes}. The proper prefixes of the twelve
orders form $13$ sets of labels, so each rectangle contributes $13$
pairs.}
\begin{tabular}{@{}lrrr@{}}
\toprule
Case & Rectangles & Interval bounds & Lemma~\ref{lem:boundary-prefixes}\\
\midrule
2 &33&427&2\\
3 &20&237&23\\
4 &21&250&23\\
\bottomrule
\end{tabular}

\label{tab:mixed-counts}
\end{table}

\begin{proposition}[One sector contains the arc]
\label{prop:mixed}
In Cases~2--4 the subcases $\sR_U\sC_V$ and $\sC_U\sR_V$ contain no
simple polygonal counterexample.
\end{proposition}

\begin{proof}
In subcase $\sR_U\sC_V$, Lemma~\ref{lem:witness-bridge}(c) gives the
twelve orders \eqref{eq:mixed-twelve}, Lemma~\ref{lem:common-mixed} the
coefficients and the value of $\Phi$, and the computation together with
Lemma~\ref{lem:boundary-prefixes} shows that every proper prefix lies in
the unit disk on the whole domain. Hence every chain has length greater
than one, contrary to Lemma~\ref{lem:weak-order-enclosure}.

For subcase $\sC_U\sR_V$, reverse the parameter and relabel
\[
 (A',B',C',D',E',F')=(F,E,D,C,B,A).
\]
In Cases~2 and~4 apply the isometry $x\mapsto R_{-\vartheta}(x-V)$, which
exchanges the apices and gives $\vartheta'=360^\circ-\vartheta$. In
Case~3 apply $x\mapsto\refl R_{-\vartheta}(x-V)$, which gives
$\vartheta'=\vartheta$ and exchanges the roles of the two caps. These isometries
preserve the incidences of each case, so the first part applies.
\end{proof}

\subsection{Conclusion of the proof}
\label{subsec:assembly}

\begin{proof}[Proof of Theorem~\ref{thm:H15}]
The set $H_{15}$ is compact and convex (Section~\ref{sec:trim}), and
Proposition~\ref{prop:area} gives its area. For the covering statement we
check the hypothesis of Proposition~\ref{prop:certificate-criterion}.
Case~1 is excluded by Proposition~\ref{prop:case1-reflection}. In
Cases~2--4, subcase $\sC_U\sC_V$ is excluded by
Propositions~\ref{prop:case2-geometric} and~\ref{prop:case34}, which give
chains of length greater than one, and the subcases $\sC_U\sR_V$ and
$\sR_U\sC_V$ by Proposition~\ref{prop:mixed}.
\end{proof}

\section{Deeper trims and minimality}
\label{sec:threshold}

By Theorem~\ref{thm:H15} and the nesting in Proposition~\ref{prop:area},
$H_t$ is a cover for every $15^\circ\le t<30^\circ$. It is natural to ask
whether one can trim further, and whether $H_{15}$ is minimal. We cannot
answer either question, but we show that $H_t$ is not a cover for
$t\le3.6^\circ$ (Proposition~\ref{prop:small-trims}).

For $t<15^\circ$ a counterexample still has witnesses in both caps, but
the open caps now overlap (Figure~\ref{fig:trim}(c)): $\pt{15}$ lies in
both, because
\[
 \cos(t/2)>\cos(15^\circ-t/2),
\]
so the two witnesses may coincide. In Lemma~\ref{lem:apex} the
corresponding parameters are $\alpha=t/2$ and
$h=\cos(15^\circ-t/2)<\cos\alpha$, so the leading term $h/\cos\alpha$ is
less than one. Lemma~\ref{lem:apex} therefore no longer gives a length
greater than one, and the estimates of Section~\ref{sec:proof} would need
new ideas below $15^\circ$ as well. The same happens for a subset of
$\PiSector$ that is symmetric in the bisector, such as a trim on both
sides: it is its own mirror image, so an arc in a unit sector with no
copy in it need have only one point outside it. Near $t=15^\circ$
each degree of trimming saves
\[
 \frac{d}{dt}\area(H_t)\Big|_{t=15^\circ}
 =\frac{\pi}{360}(1-\cos15^\circ)
 \approx2.974\times10^{-4}
\]
in area (Proposition~\ref{prop:area}).

\subsection{Minimality}
\label{subsec:minimality}

We call a compact convex cover \emph{inclusion-minimal} if no proper
compact convex subset of it is a cover.

\begin{proposition}[Unit segments in $H_{15}$]
\label{prop:minimality-apex}
Every subset of $H_{15}$ that contains a unit segment contains $O$ and at
least one point of
\[
 \{\pt{\theta}:0^\circ\le\theta\le15^\circ\}\cup\{\pt{30}\}.
\]
Moreover, no homothetic copy of $H_{15}$ with ratio less than one is a
cover.
\end{proposition}

\begin{proof}
Let $x=r_1\pt{\theta_1}$ and $y=r_2\pt{\theta_2}$ belong to $\PiSector$,
with $1\ge r_1\ge r_2>0$. Since $|\theta_1-\theta_2|\le30^\circ$,
\[
 |x-y|^2
 =r_1^2+r_2^2-2r_1r_2\cos(\theta_1-\theta_2)
 \le r_1^2+r_2(r_2-\sqrt3r_1)<r_1^2\le1.
\]
Hence, if two points of $\PiSector$ are at distance one, one of them is
$O$ and the other has norm one; the points of $H_{15}$ of norm one are those listed. Apply
this to the endpoints of a unit segment. Since $H_{15}$ has diameter one,
a homothetic copy with ratio less than one has smaller diameter and
contains no unit segment.
\end{proof}

\begin{proposition}[Cuts at extreme points]
\label{prop:minimality-cuts}
The cover $H_{15}$ is inclusion-minimal if and only if every member of
the two families
\[
 H_{15}\cap\{x:\ip{x}{\pt{\theta}}\le1-\eta\},
 \qquad 0^\circ\le\theta\le15^\circ,\quad \eta>0,
\]
\[
 H_{15}\cap\{x:\ip{x}{\pt{90}}\le\tfrac12-\eta\},
 \qquad \eta>0,
\]
fails to cover some unit arc.
\end{proposition}

\begin{proof}
A member of either family that is a cover is a proper compact convex
subcover. Conversely, let $K\subsetneq H_{15}$ be a compact convex cover.
Since
\[
 H_{15}=\conv\bigl(\{O,\pt{30}\}
          \cup\{\pt{\theta}:0^\circ\le\theta\le15^\circ\}\bigr),
\]
$K$ omits one of these extreme points, and by
Proposition~\ref{prop:minimality-apex} the omitted point is not $O$. If $K$ omits
$\pt{\theta}$, then compactness, together with the fact that $\pt{\theta}$
is the only point of $H_{15}$ with $\ip{x}{\pt{\theta}}=1$, gives
$\max_{x\in K}\ip{x}{\pt{\theta}}<1$, so $K$ lies in a member of the first
family. If $K$ omits $\pt{30}$, the same argument with the functional
$\ip{x}{\pt{90}}$, which $\pt{30}$ alone maximizes on $H_{15}$, places $K$
in a member of the second family. In either case that member contains
the cover $K$ and is therefore a cover.
\end{proof}

\begin{openproblem}
\label{op:inclusion-minimality}
Is $H_{15}$ inclusion-minimal among compact convex covers of unit arcs?
Equivalently, does every set in the two families of
Proposition~\ref{prop:minimality-cuts} fail to contain a congruent copy of
some unit arc?
\end{openproblem}

We consider only compact subsets because the reduction to polygonal arcs
\cite[Corollary~5]{wetzel-wichiramala2010} requires compactness.

Our method does not apply as it stands to the first family. A set $K$ of
that family omits points $\pt\phi$ with $0^\circ<\phi<15^\circ$; these
lie in the open cap $C_{15}^-$, hence outside both $K$ and $\refl(K)$, so
one point can serve as both witnesses, as for $t<15^\circ$. The set
conjectured in \cite{movshovich2025}, placed as in
Section~\ref{sec:intro}, omits $\pt{30}$ and every $\pt\theta$ with
$\theta>0^\circ$, so if it lies in $H_{15}$, it lies in members of both
families.

The unit segment, a certain $z$-arc $\arcZ$ and the
$60^\circ$ staples of \cite{movshovich-wetzel2017} all have congruent
copies in a proper subset $K_0$ of $H_{15}$ that is not a cover
(Example~\ref{ex:tests}), so they do not decide Open
Problem~\ref{op:inclusion-minimality}. Nor can
finitely many test arcs decide in general whether a convex set is a
cover: for each $n$ there is a convex set that contains a congruent copy
of every unit arc made of $n$ segments but is not a cover
\cite{panraksa-wetzel-wichiramala2007}.

\begin{figure}[tbp]
\centering
\includegraphics[width=0.78\linewidth]{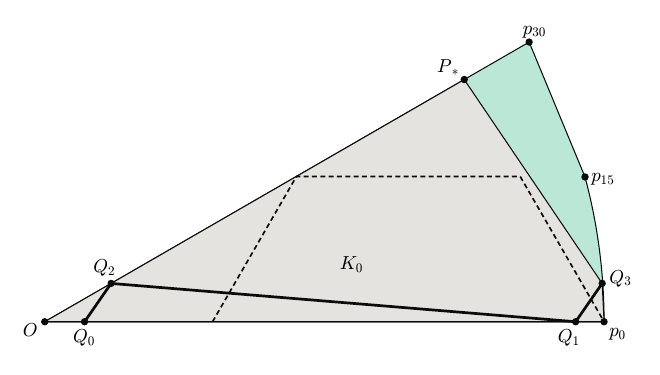}
\caption{Example~\ref{ex:tests}. The quadrilateral
$K_0=\conv\{O,\pt0,Q_3,P_*\}$ (gray) inside $H_{15}$; the part of $H_{15}$
outside $K_0$ is green. The set $K_0$ contains the copy $Q_0Q_2Q_1Q_3$ of
$\arcZ$ (thick) and the staple for $d=0.7$, with vertices $(0.3,0)$,
$0.6P_*$, $0.6P_*+0.4\,\pt0$ and $\pt0$ (dashed).}
\label{fig:minimality-tests}
\end{figure}

\begin{example}
\label{ex:tests}
Let $\arcZ$ be the unit $z$-arc $Z_0Z_1Z_2Z_3$, whose first and last segments
are parallel and point the same way, with
\[
  Z_0=(0,0),\qquad
  Z_1=\Bigl(\tfrac1{12},0\Bigr),\qquad
  Z_2=\Bigl(\tfrac12,\tfrac{5\sqrt3}{12}\Bigr),\qquad
  Z_3=\Bigl(\tfrac7{12},\tfrac{5\sqrt3}{12}\Bigr).
\]
Movshovich \cite{movshovich2021} studies such arcs in Wetzel's sector. The
convex hull of $\arcZ$ is the parallelogram spanned by $Z_1-Z_0$ and
$Z_2-Z_0$, with
sides $1/12$ and $\sqrt{111}/12$ and angle $\beta$, where
$\cos\beta=6/\sqrt{111}$ and $\sin\beta=5\sqrt3/\sqrt{111}$, so that
$\beta=55.28\ldots^\circ$. Write $\bar z$ for the reflection of $z$ in the
first coordinate axis, and put
\[
 Q_0=\Bigl(\frac{3}{4\sqrt{111}},0\Bigr).
\]
The map $z\mapsto R_\beta\bar z$, with matrix
\[
 \frac1{\sqrt{111}}
 \begin{pmatrix}6&5\sqrt3\\5\sqrt3&-6\end{pmatrix},
\]
is the reflection in the line through $O$ at angle $\beta/2$, and it
exchanges the directions of $Z_1-Z_0$ and $Z_2-Z_0$. The map
$z\mapsto Q_0+R_\beta\bar z$ sends $Z_0,Z_1,Z_2,Z_3$ to
$Q_0,Q_2,Q_1,Q_3$, where
\begin{align*}
Q_1&=\Bigl(\frac{10}{\sqrt{111}},0\Bigr),
&
Q_2&=\Bigl(\frac{5}{4\sqrt{111}},
\frac{5\sqrt3}{12\sqrt{111}}\Bigr),
&
Q_3&=\Bigl(\frac{21}{2\sqrt{111}},
\frac{5\sqrt3}{12\sqrt{111}}\Bigr).
\end{align*}
The points $Q_0$ and $Q_1$ lie on the radius $[O,\pt0]$. The point $Q_2$
lies on the radius $[O,\pt{30}]$, since $y(Q_2)/x(Q_2)=1/\sqrt3$. Finally
\[
  |Q_3|^2=\frac{5317}{5328}<1,
  \qquad
  \frac{y(Q_3)}{x(Q_3)}
  =
  \frac{5\sqrt3}{126}
  <
  2-\sqrt3
  =
  \tan15^\circ,
\]
so $Q_3\in\Sigma_{15}$. By convexity the image of the hull, and with it
this copy of $\arcZ$, lies in $H_{15}$.

Let $P_*=(3/4,\sqrt3/4)$ be the foot of the perpendicular from $\pt0$ to
the radius $[O,\pt{30}]$, so that $|\pt0-P_*|=1/2$, and put
\[
 K_0=\conv\{O,\pt0,Q_3,P_*\}\subseteq H_{15}
\]
(Figure~\ref{fig:minimality-tests}). The set $K_0$ contains the unit
segment $[O,\pt0]$ and the copy of $\arcZ$ above, since $Q_0,Q_1\in[O,\pt0]$
and $Q_2=5P_*/(3\sqrt{111})\in[O,P_*]$. It also contains a copy of every
unit $60^\circ$ staple of \cite{movshovich-wetzel2017}. There is one for
each $d$ with
$1/2\le d<1$: it consists of three sides, of lengths $1-d$, $2d-1$ and
$1-d$, of an isosceles trapezoid with base angles $60^\circ$ whose fourth
side has length $d$. The vertices
\[
 (1-d,0),\quad 2(1-d)P_*,\quad 2(1-d)P_*+(2d-1)\pt0,\quad \pt0
\]
give such a staple, and all four lie in $\conv\{O,P_*,\pt0\}$. Of the
vertices of $K_0$ only $\pt0$ has norm one,
so by strict convexity of the disk $K_0$ meets the unit circle only at
$\pt0$, and $K_0\ne H_{15}$. The shoelace formula gives
\[
 \area(K_0)=\frac{\sqrt3\,x(Q_3)+y(Q_3)}8
 =\frac{131}{96\sqrt{37}}=0.2243361\ldots,
\]
which is below the lower bound of \cite{khandhawit-pagonakis-sriswasdi2013},
so $K_0$ is not a cover. The boundary of $K_0$ consists of the radius
$[O,\pt0]$, the segment $[O,P_*]$ and the convex arc $\pt0Q_3P_*$, so it
has the form conjectured in \cite{movshovich2025} for a convex cover of
least area, with the conjectured distance $1/2$ from $\pt0$ to $P_*$.
\end{example}

\subsection{Trims that are not covers}
\label{subsec:noncovers}

The arc $\arcZ$ of Example~\ref{ex:tests} shows that $H_t$ is not a cover
when $t$ is small. A congruent copy of $\arcZ$ is its image under a map
$z\mapsto x_0+R_\theta z$ or $z\mapsto x_0+R_\theta\bar z$, with $\bar z$
as in Example~\ref{ex:tests}; we call the copy \emph{direct} in the first
case and \emph{reflected} in the second.
The sides of its convex hull have directions $\theta$ and
$\theta+\beta$ for a direct copy, and $\theta$ and $\theta-\beta$ for a
reflected one. In the formulas below the upper sign refers to direct
copies and the lower sign to reflected ones.

\begin{lemma}[A test with the arc $\arcZ$]
\label{lem:ztest}
Let $0^\circ\le\chi\le30^\circ$ and $h\in\R$, and put
\[
 f_\chi(\phi)=\sin(120^\circ-\chi)\,|\sin\phi|
 +\cos\chi\,|\cos(\phi-120^\circ)|+\tfrac12|\cos(\phi-\chi)|.
\]
If the triangle
\[
 \{x:\ \ip{x}{\pt{270}}\le0,\ \ip{x}{\pt{120}}\le0,\
 \ip{x}{\pt\chi}\le h\}
\]
contains a copy of $\arcZ$ given by $\theta$ as above, then
\begin{equation}
 \frac1{12}\,f_\chi(\theta)
 +\frac{\sqrt{111}}{12}\,f_\chi(\theta\pm\beta)\le h.
 \label{eq:ztest}
\end{equation}
For either sign the left-hand side has period $180^\circ$ in $\theta$,
and it is concave on every interval whose interior contains none of the
points
\begin{equation}
 0^\circ,\quad 30^\circ,\quad 90^\circ+\chi,\quad
 \mp\beta,\quad 30^\circ\mp\beta,\quad
 90^\circ+\chi\mp\beta
 \pmod{180^\circ}.
 \label{eq:zbreak}
\end{equation}
\end{lemma}

\begin{proof}
The triangle is convex, so it contains the convex hull of the copy, a
parallelogram whose sides have lengths $1/12$ and $\sqrt{111}/12$ and
directions $\theta$ and $\theta\pm\beta$. Let $x_1$ be its center. For a
unit vector $n$, the maximum of $\ip{x}{n}$ over the parallelogram is
\[
 \ip{x_1}{n}+\frac1{24}\bigl|\ip{\pt\theta}{n}\bigr|
 +\frac{\sqrt{111}}{24}\bigl|\ip{\pt{\theta\pm\beta}}{n}\bigr| .
\]
For $n=\pt{270}$, $\pt{120}$ and $\pt\chi$ this is at most $0$, $0$ and
$h$. Multiply these three inequalities by $\sin(120^\circ-\chi)$,
$\cos\chi$ and $1/2$, which are positive, and add. The terms in $x_1$
cancel, because
\[
 \sin(120^\circ-\chi)\,\pt{270}+\cos\chi\,\pt{120}+\tfrac12\,\pt\chi=0,
\]
and since
$|\ip{\pt\phi}{\pt{270}}|=|\sin\phi|$,
$|\ip{\pt\phi}{\pt{120}}|=|\cos(\phi-120^\circ)|$ and
$|\ip{\pt\phi}{\pt\chi}|=|\cos(\phi-\chi)|$, the sum is half of
\eqref{eq:ztest}. Each of these three functions of $\phi$ has period
$180^\circ$ and is concave between consecutive zeros, which are at
$0^\circ$, $30^\circ$ and $90^\circ+\chi$ modulo $180^\circ$. Hence
$f_\chi$ is concave on every interval whose interior contains none of
these three points, and the last statement follows.
\end{proof}

\begin{proposition}
\label{prop:small-trims}
The arc $\arcZ$ has no congruent copy in $H_t$ for $0^\circ<t\le3.6^\circ$,
and it has one for $3.7^\circ\le t<30^\circ$. In particular, $H_t$ is not
a cover for $t\le3.6^\circ$, and neither is the triangle
$\conv\{O,\pt0,\pt{30}\}$.
\end{proposition}

\begin{proof}
By the nesting in Proposition~\ref{prop:area} it suffices to prove the
first statement for $t=3.6^\circ$ and the second for $t=3.7^\circ$. For
$t=3.6^\circ$ we have $m_t=16.8^\circ$ and $\delta_t=13.2^\circ$, so $H_t$
lies in the triangle of Lemma~\ref{lem:ztest} with $\chi=16.8^\circ$ and
$h=\cos13.2^\circ=0.97357\ldots$. By the lemma, the left-hand side of
\eqref{eq:ztest} is smallest at one of the points \eqref{eq:zbreak}.
Table~\ref{tab:ztest} lists its values there. They all exceed $h$, so
\eqref{eq:ztest} fails for every $\theta$ and either sign, and $\arcZ$ has
no copy in $H_t$.

For $t=3.7^\circ$ we have $m_t=16.85^\circ$ and $\delta_t=13.15^\circ$.
The copy $Q_0Q_2Q_1Q_3$ of Example~\ref{ex:tests} has $Q_0$, $Q_1$ and
$Q_2$ on the two radii and $Q_3\in\PiSector$, with
\[
 \ip{Q_3}{\pt{16.85}}=0.97368\ldots<\cos13.15^\circ=0.97377\ldots,
\]
so it lies in $H_t$. Finally, the triangle is contained in every $H_t$,
by \eqref{eq:triangular-core}.
\end{proof}

\begin{table}[htbp]
\centering
\caption{The left-hand side of \eqref{eq:ztest} for $\chi=16.8^\circ$ at
the points \eqref{eq:zbreak}, for direct and for reflected copies,
rounded down to five decimal places; here
$h=\cos13.2^\circ=0.97357\ldots$.}
\begin{tabular}{@{}lcccccc@{}}
\toprule
Direct: $\theta$ & $0^\circ$ & $30^\circ$ & $106.8^\circ$ &
$180^\circ-\beta$ & $210^\circ-\beta$ & $106.8^\circ-\beta$\\
value & $1.48501$ & $1.78489$ & $1.40289$ & $0.99951$ & $0.98592$ &
$1.76359$\\
\midrule
Reflected: $\theta$ & $0^\circ$ & $30^\circ$ & $106.8^\circ$ &
$\beta$ & $30^\circ+\beta$ & $106.8^\circ+\beta$\\
value & $1.75508$ & $1.46290$ & $1.49352$ & $0.97387$ & $1.01648$ &
$1.75499$\\
\bottomrule
\end{tabular}

\label{tab:ztest}
\end{table}

The smallest entry of Table~\ref{tab:ztest} belongs to the reflected copy
with $\theta=\beta$, whose long side lies along the radius $[O,\pt0]$;
this is the position of the copy in Example~\ref{ex:tests}. A program in
the repository checks the inequalities of this proof, and those of
Remark~\ref{rem:pentagon} below, in exact arithmetic.

\begin{openproblem}
\label{op:tstar}
Determine
\[
 t^*=\inf\{t\in(0^\circ,30^\circ):H_t\text{ covers every unit arc}\}.
\]
Does some $H_t$ with $t<15^\circ$ cover all unit arcs?
\end{openproblem}

Theorem~\ref{thm:H15} gives $t^*\le15^\circ$, and
Proposition~\ref{prop:small-trims} gives $t^*\ge3.6^\circ$. Since $\arcZ$ has
a copy in $H_t$ for $t\ge3.7^\circ$, this arc cannot give a lower bound
above $3.7^\circ$.

\begin{remark}
\label{rem:pentagon}
A subset of the sector with larger area than $H_{15}$ need not be a
cover. The pentagon $\conv\{O,\pt0,\pt{10},\pt{20},\pt{30}\}$, in which
three chords replace the arc, has area
$\frac32\sin10^\circ=0.2604722\ldots$, more than that of $H_{15}$, but it
contains no copy of $\arcZ$. Since the pentagon is symmetric in the bisector
of $\PiSector$, it suffices to consider direct copies. It lies in the
triangle of
Lemma~\ref{lem:ztest} for $(\chi,h)=(15^\circ,\cos5^\circ)$ and for
$(25^\circ,\cos5^\circ)$. For direct copies, the left-hand side of
\eqref{eq:ztest} with $\chi=15^\circ$ is at least $1.00849$ on
$[0^\circ,140^\circ]$, its smallest value there being at
$\theta=180^\circ-\beta$, and with $\chi=25^\circ$ it is at least
$0.99879$ on $[140^\circ,180^\circ]$, its smallest value being at
$\theta=210^\circ-\beta$; by the concavity in Lemma~\ref{lem:ztest} one
need only compare the values at the points \eqref{eq:zbreak} in these
intervals and at their endpoints. Both bounds exceed
$\cos5^\circ=0.99619\ldots$.
\end{remark}

\subsection{Covers of diameter one}
\label{subsec:diameter}

Every cover has diameter at least one, and by the proof of
Proposition~\ref{prop:minimality-apex} no two points of $\PiSector$ are
farther apart than one. So a subset of $\PiSector$ that is a cover has
diameter one. The sets considered in this section and the set
conjectured in \cite{movshovich2025} are all subsets of $\PiSector$.

Among compact convex covers of diameter one there is one of least area.
To see this, translate such covers so that they contain $O$; they then
lie in the unit disk. By the Blaschke selection theorem, a sequence of
them whose areas tend to the infimum has a subsequence that converges in
the Hausdorff metric to a compact convex set $K$. Area and diameter are
continuous under this convergence. Moreover $K$ is a cover: the
isometries that carry a given unit arc into the members of the
subsequence have a limit point, and this isometry carries the arc into
$K$. For convex covers with no condition on the diameter, the existence
of one of least area is stated in \cite{norwood-poole-laidacker1992}.

\begin{openproblem}
\label{op:diameter}
What is the least area of a convex cover of diameter one? Is it less
than $\area(H_{15})$? Is it equal to the least area of a convex cover?
\end{openproblem}

By Theorem~\ref{thm:H15} and
\cite{khandhawit-pagonakis-sriswasdi2013}, this least area lies between
$0.232239$ and $\area(H_{15})$. It is less than $\area(H_{15})$ if
$H_{15}$ is not inclusion-minimal, and in particular if
$t^*<15^\circ$. If the conjecture of \cite{movshovich2025} holds, then
the answer to the last question is yes. The survey poses the same
problem with the least possible width in place of the least possible
diameter \cite[Problem~2]{movshovich2025}.

\section*{Disclosure statement}
The authors report no competing interests.

\section*{Funding}
No funding was received for this research.

\section*{Data and code availability}
The coefficient data and the Python programs used in the proof, with
instructions and tests, are available at
\begin{center}
\url{https://github.com/chatchawanpan-dev/wetzel-sector-trim-computations}
\end{center}
and a copy is included with the arXiv version of this paper as ancillary
files. The programs that check the proof use only the Python standard
library. From the root folder of the repository, the command
\texttt{python3 verification/verify\_proof\_essential.py} checks the
computations in the proof of Theorem~\ref{thm:H15}, and
\texttt{python3 verification/Section6/verify\_zarc\_test.py} checks those
of Section~\ref{subsec:noncovers}. The file \texttt{REPRODUCIBILITY.md}
lists the other commands, the parameter domains and the arithmetic
conventions. A second program, described in \texttt{README.md}, was
written separately from the statements of the paper. It repeats the interval
checks of Section~\ref{subsec:balanced-normals} with the interval
arithmetic of the library mpmath, and it arrives at the counts of
Table~\ref{tab:mixed-counts}; it needs NumPy and mpmath. The folder \texttt{exploration} holds the two programs for
the experiments described in Section~\ref{subsec:method}; they need
NumPy, SciPy and CVXPY. The figures are
drawn by \texttt{make\_figures.py}, also in the repository, which needs
NumPy and Matplotlib.

\section*{Declaration of generative AI use}
In preparing this paper the authors used ChatGPT~6.1 Sol (OpenAI) and
Claude Fable~5.2 (Anthropic). The tools were used to help find and check
arguments, to draft and revise parts of the text, among them parts of
Sections~\ref{sec:intro}, \ref{sec:reduction} and~\ref{sec:threshold} and the
figure captions, and to write and test the programs in the repository,
including the program that draws the figures. The authors have checked
the mathematics, the programs and the text, and they take full
responsibility for the content of the paper.

\appendix

\section{Prefix bounds when both sectors contain the arc}
\label{sec:capcap-circles}

We use $e_0$, $e_1$ and $\kappa$ from Section~\ref{sec:proof}, and recall
that
\[
 b=2\cos15^\circ-1=\frac{\sqrt6+\sqrt2}{2}-1.
\]
Then
\begin{equation}
 \kappa\pt{22.5}=e_0+b e_1,
 \qquad \kappa\pt{7.5}=b e_0+e_1.
 \label{eq:normalized-cap-b}
\end{equation}

\subsection{Case 2}
\label{subsec:case2-global-bounds}

Use the coefficients of Proposition~\ref{prop:case2-geometric} and put
$x=\sin\vartheta$ and $y=-\cos\vartheta$, so that $-1/2\le x\le1/2$,
$y\ge0$ and $x^2+y^2=1$. Moreover, $y\le1-x^2/2$, because the
right-hand side is positive and its square exceeds $1-x^2$ by $x^4/4$.

The four orders have nine distinct proper prefixes. Four of them have
norms $|a|$, $|Q|$, $|1+c|$ and $|Q|$, all at most $5/6$. The other five are
\[
\begin{array}{c|c}
\{A,X_V\}&a e_0+R_\vartheta Q\\
\{A,B,X_V\}&a e_0-b e_1+R_\vartheta Q\\
\{A,B,C,X_V\}&a e_0-b e_1+R_\vartheta(Q+c e_0)\\
\{A,\ldots,D,X_V\}&-Q+R_\vartheta(Q+c e_0)\\
\{A,\ldots,E,X_V\}&-Q+(1+c)R_\vartheta e_0.
\end{array}
\]
Denote their squared norms by $N_1,\ldots,N_5$. Since $a(-x)=-1-c(x)$ and
$c(-x)=-1-a(x)$, rotation through $-\vartheta$ and negation give
$N_5(x,y)=N_1(-x,y)$ and $N_4(x,y)=N_2(-x,y)$, so it suffices to bound
the first three.

Write
\[
 Q=(q_1,q_2),\qquad q_1=b/2,\qquad q_2=1-\sqrt3\,q_1,\qquad
 k=q_1^2+(q_2-1/3)^2.
\]
The enclosures
\[
 1.73205<\sqrt3<1.73206,\qquad 0.9318<b<0.9319,
\]
which follow by squaring and from $(b+1)^2=2+\sqrt3$ with $b>0$, imply
$q_1,q_2>0$, $q_1^2+q_2^2-q_2/2>0$, $k>1/36$, and
\[
 4k<19/20,\qquad |Q|^2=q_1^2+q_2^2<13/50.
\]

For the first two norms, substitute $y\le1-x^2/2$ in their expansions.
This gives upper bounds because the coefficients of $y$,
\[
 \frac{2q_2(2+x)}{3},\qquad
 2\left(q_1^2+q_2^2-\frac{q_2}3+\frac{q_2x}{3}\right),
\]
are positive on $[-1/2,1/2]$. The resulting cubic has coefficient
$-q_2/3$ at $x^3$, and $-x^3\le x^2/2$ gives
\[
 N_i\le A_i x^2+B_i x+C_i\qquad (i=1,2),
\]
where
\[
\begin{array}{c|c|c|c}
i&A_i&B_i&C_i\\\hline
1&\frac19-\frac{q_2}2-\frac{2q_1}3
 &\frac{2q_2}3-\frac{4q_1}3+\frac49
 &q_1^2+(q_2+\frac23)^2\\[2pt]
2&-q_2^2+\frac{q_2}2-q_1^2-\frac{2q_1}3+\frac19
 &\frac{4q_2}3+\frac{2q_1}3-\frac29
 &4q_2^2-\frac{4q_2}3+4q_1^2+\frac19.
\end{array}
\]
The same enclosures give
\[
\begin{array}{c|ccc}
i&-A_i\ge\hat A_i&|B_i|\le\hat B_i&C_i\le\hat C_i\\\hline
1&29/100&1/20&957/1000\\
2&7/20&7/20&109/125,
\end{array}
\]
and maximizing $-\hat A_i x^2+\hat B_i|x|+\hat C_i$ yields
\[
 N_1\le\frac{957}{1000}+\frac1{464}
       =\frac{55631}{58000}
       <\left(\frac{49}{50}\right)^2,
 \qquad
 N_2\le\frac{109}{125}+\frac7{80}
       =\frac{1919}{2000}
       <\left(\frac{49}{50}\right)^2.
\]
For the third prefix a direct expansion gives
\[
 N_3=2k(1+y)+\frac{2x^2}{9}(1-y)-\frac{4q_1 x^2}{3}.
\]
As $k>x^2/9$ and $y\le1$,
\[
 N_3\le4k-\frac{4q_1 x^2}{3}\le4k
 <\frac{19}{20}<\left(\frac{49}{50}\right)^2.
\]
Together with the constant prefixes and the symmetry, these
estimates give the bound used in Proposition~\ref{prop:case2-geometric}.
A program in the repository checks the polynomial identities, modulo
$x^2+y^2=1$, and the rational inequalities above in exact arithmetic.

\subsection{Cases 3 and 4}
\label{subsec:case34-circles}

Let $(\phi_U,\phi_V)$ be as in Proposition~\ref{prop:case34}, write
$h_0(\psi)e_0+h_1(\psi)e_1=-\kappa\pt\psi$, and take $\varepsilon=1$ in
Case~3 and $\varepsilon=0$ in Case~4. The coefficients are
\[
\begin{array}{c|cc}
 & [150^\circ,180^\circ]&[180^\circ,210^\circ]\\\hline
\text{Case 3: }(a,c)&(-13/20,-5/16)&(-1/2,-1/2)\\
\text{Case 4: }(a,c)&(-3/4,-1/2)&(-1/2,-1/4);
\end{array}
\]
at $\vartheta=180^\circ$ either column may be used, together with its
pair $(\phi_U,\phi_V)$. They are assigned as
\[
\begin{array}{c|c@{\quad}c|c}
A&a e_0&D&cR_\vartheta e_\varepsilon\\
B&h_1(\phi_U)e_1&E&h_{1-\varepsilon}(\phi_V)R_\vartheta e_{1-\varepsilon}\\
C&(h_0(\phi_U)-a)e_0&F&(h_\varepsilon(\phi_V)-c)R_\vartheta e_\varepsilon\\
X_U&\kappa\pt{\phi_U}&X_V&\kappa\pt{\vartheta+\phi_V}.
\end{array}
\]

The four orders have nine proper prefixes. Four of them have norms
$|a|$, $\kappa$, $|h_\varepsilon(\phi_V)-c|$ and $\kappa$, all at most
$3/4$. The other five are $u+R_\vartheta v$ with
\[
\begin{array}{c|cc}
\text{prefix}&u&v\\\hline
\{A,X_V\}&a e_0&\kappa\pt{\phi_V}\\
\{A,B,X_V\}&a e_0+h_1(\phi_U)e_1&\kappa\pt{\phi_V}\\
\{A,B,C,X_V\}&-\kappa\pt{\phi_U}&\kappa\pt{\phi_V}\\
\{A,\ldots,D,X_V\}&-\kappa\pt{\phi_U}&\kappa\pt{\phi_V}+c e_\varepsilon\\
\{A,\ldots,E,X_V\}&-\kappa\pt{\phi_U}&-(h_\varepsilon(\phi_V)-c)e_\varepsilon.
\end{array}
\]
For the middle row, the directed angle $\vartheta+\phi_V-\phi_U$ lies in
$[135^\circ,165^\circ]$ on the first interval and in
$[195^\circ,225^\circ]$ on the second, so the angle between the two
vectors is at most $165^\circ$, and
\[
 \bigl|\kappa\pt{\vartheta+\phi_V}-\kappa\pt{\phi_U}\bigr|
 \le2\kappa\sin82.5^\circ=1,
\]
with equality only at $\vartheta=180^\circ$. For the other four rows the
criterion in the proof of Proposition~\ref{prop:case34} locates the
maximum of the norm, and Table~\ref{tab:circle-maxima} lists the sixteen
maxima. Each is less than $M=999/1000$. These are finite inequalities in
$b$ and $\sqrt3$, proved with the enclosures of
Appendix~\ref{subsec:case2-global-bounds}; for an interior maximum the
bound $|u|+|v|<M$ follows from
\[
 M^2-|u|^2-|v|^2>0,
 \qquad \bigl(M^2-|u|^2-|v|^2\bigr)^2-4|u|^2|v|^2>0.
\]
A program in the repository checks these inequalities in exact
arithmetic. Hence all nine prefixes lie in the unit disk.

\begin{table}[!htb]
\centering
\caption{The largest norm of each of the four prefixes $u+R_\vartheta v$ of
Appendix~\ref{subsec:case34-circles} that move on circles, rounded up to
four decimal places, and the
value of $\vartheta$ at which it occurs, rounded to $0.1^\circ$; an
asterisk marks an interior maximum, of norm $|u|+|v|$.}
\begin{tabular}{@{}llcccc@{}}
\toprule
Case & $\vartheta$ & $\{A,X_V\}$ & $\{A,B,X_V\}$ & $\{A,\ldots,D,X_V\}$ & $\{A,\ldots,E,X_V\}$\\
\midrule
3 & $[150^\circ,180^\circ]$ & $0.8732$ ($180^\circ$) & $0.9960^*$ ($153.9^\circ$) & $0.9854^*$ ($158.1^\circ$) & $0.9042$ ($180^\circ$)\\
3 & $[180^\circ,210^\circ]$ & $0.9511$ ($210^\circ$) & $0.9814$ ($180^\circ$) & $0.9956$ ($180^\circ$) & $0.8868$ ($210^\circ$)\\
4 & $[150^\circ,180^\circ]$ & $0.9569$ ($180^\circ$) & $0.9738^*$ ($165.5^\circ$) & $0.9956$ ($180^\circ$) & $0.8868$ ($150^\circ$)\\
4 & $[180^\circ,210^\circ]$ & $0.9511$ ($210^\circ$) & $0.9814$ ($180^\circ$) & $0.9738^*$ ($194.5^\circ$) & $0.9569$ ($180^\circ$)\\
\bottomrule
\end{tabular}

\label{tab:circle-maxima}
\end{table}

\end{document}